\documentclass[12pt,a4paper,reqno]{article}%
\usepackage{amssymb}
\usepackage{amsmath}
\usepackage{amsfonts}
\usepackage{graphicx}
\usepackage{mathrsfs}
\usepackage[all]{xy}
\usepackage{hyperref}%
\providecommand{\U}[1]{\protect\rule{.1in}{.1in}}
\SelectTips{cm}{}
\newtheorem{theorem}{Theorem}
\newtheorem{lemma}[theorem]{Lemma}

\newtheorem{definition}{Definition}
\newtheorem{example}{Example}
\newtheorem{proposition}[theorem]{Proposition}
\newtheorem{corollary}[theorem]{Corollary}

\newtheorem{problem}{Problem}
\newenvironment{proof}[1][Proof]{\noindent\textbf{#1.} }{\ \rule{0.5em}{0.5em}}
\begin{document}

\title{\textbf{The Hamilton-Jacobi Formalism for Higher Order Field Theories}}
\author{\textsc{{L.~Vitagliano}\thanks{\textbf{e}-\textit{mail}:
\texttt{lvitagliano@unisa.it}}}\\{\small {DMI, Universit\`a degli Studi di Salerno, and}}\\{\small {Istituto Nazionale di Fisica Nucleare, GC Salerno}}\\{\small {Via Ponte don Melillo, 84084 Fisciano (SA), Italy}}}
\maketitle

\begin{abstract}
We extend the geometric Hamilton-Jacobi formalism for hamiltonian mechanics to
higher order field theories with regular lagrangian density. We also
investigate the dependence of the formalism on the lagrangian density in the
class of those yelding the same Euler-Lagrange equations.

\end{abstract}

\section*{Introduction}

The Hamilton-Jacobi (HJ) formalism is a cornerstone of the calculus of
variations (see for instance \cite{r73}) and the theory of Hamiltonian
systems. Moreover, it is a first, important step through the quantization of a
mechanical system (see, for instance \cite{ems04}, see also \cite{mmm09}). HJ
formalism can be readily extended to first order Lagrangian (and Hamiltonian)
field theories \cite{r73,h35}. Moreover, both its original version and its
first order, field theoretic extension posses an effective geometric
formulation in terms of symplectic \cite{vk81} and multisymplectic
\cite{pr02,pr02b,dL...08} geometry respectively. Finally, in \cite{c...06} the
authors formulate in geometric terms a generalized HJ problem depending on the
sole equations of motions (and not directly on the Lagrangian, nor the Hamiltonian
function itself). In particular, such generalized problem can be stated for
any SODE on the tangent bundle and any vector field on the cotangent bundle of
a configuration manifold, and thus it has a wide range of applicability. The
aim of the present paper is to formulate in geometric terms a (generalized) HJ
problem for higher order Lagrangian field theories, in view of its application
to both variational calculus and theoretical physics. Recall that higher order
Lagrangian field theory has got a very elegant geometric
formulation (see, for instance, \cite{v84}). Moreover, the Hamiltonian
formulation of Lagrangian mechanics has been recently extended to higher
order field theories (on fiber bundles) by the author \cite{v09}. Unlike
previously proposed ones, the Hamiltonian formalism of \cite{v09} is free from
any relevant ambiguity and does only depend on the action functional and the order of a Lagrangian
density up to isomorphisms. Such theory is the starting point
of the present work. It should be finally mentioned that we restrict to regular Lagrangian field theories. The HJ formalism for  general Lagrangian
field theories will be analyzed elsewhere (see \cite{c...08} for the case of a
singular mechanical system).

The paper is divided into eight sections and one appendix. The first section
summarizes the notation and convention adopted throughout the paper. Section
\ref{SecConn} summarizes very well known facts about Ehresmann connections
whose geometry plays an important role in the whole paper. Section
\ref{SecConn} also contains some less standard definitions (and results) about
\emph{relative} connections. In Section \ref{SecLHForm} we present a finite
dimensional version of the formalism in \cite{v09}. Namely, Section
\ref{SecLHForm} contains some results of \cite{v09} but they are here derived
in a new, original fashion. Section \ref{SecLHForm} also contains
some original results which are presented here for the first time. In Section
\ref{SecInvLeg} we show that although the concept of \emph{Legendre transform}
has no natural generalization to higher order Lagrangian field theory, it is
still possible to give a natural, geometric definition of \textquotedblleft%
\emph{inverse Legendre transform}\textquotedblright\ for regular theories. Such
inverse Legendre transform plays a prominent role in the HJ formalism. In
Section \ref{SecGenHJForm} we extend to higher order, regular, Lagrangian
field theories (in the Hamiltonian picture) the generalized HJ theory of
\cite{c...06}. In particular we state a generalized HJ problem and
characterize its solutions. In Section \ref{SecHJForm} we state the
(non-generalized) HJ problem for higher order field theories and present
coordinate expressions of all involved geometric objects. In Section
\ref{SecEquivLag} we study the dependence of the HJ problem on the choice of a
Lagrangian in the class of those yielding the same Euler-Lagrange equations. In
particular, we find that the HJ problems determined by equivalent Lagrangians (of the same order)
are equivalent as well. Thus the HJ formalism depends on the sole action
functional (and the order of a Lagrangian density) up to isomorphisms. In Section \ref{SecExamp}, we illustrate the
higher order, field theoretic HJ problem via a specific, simple example: the
biharmonic equation. The paper in completed by an appendix in which we
characterize Hamiltonian field theories coming from (hyperregular) Lagrangian
field theories.

\section{Notations and Conventions}

In this section we collect notation and convention about some general
constructions in differential geometry that will be used in the following.

Let $N$ be a smooth manifold. If $L\subset N$ is a submanifold, we denote by
$i_{L}:L\hookrightarrow N$ the inclusion. We denote by $C^{\infty}(N)$ the
$\mathbb{R}$--algebra of smooth, $\mathbb{R}$--valued functions on $N$. We
will always understand a vector field $X$ on $N$ as a derivation $X:C^{\infty
}(N)\longrightarrow C^{\infty}(N)$. We denote by $D(N)$ the $C^{\infty}%
(N)$--module of vector fields over $N$, by $\Lambda(M)=\bigoplus_{k}%
\Lambda^{k}(N)$ the graded $\mathbb{R}$--algebra of differential forms over
$N$ and by $d:\Lambda(N)\longrightarrow\Lambda(N)$ the de Rham differential.
If $F:N_{1}\longrightarrow N$ is a smooth map of manifolds, we denote by
$F^{\ast}:\Lambda(N)\longrightarrow\Lambda(N_{1})$ the pull-back via $F$. We will
understand everywhere the wedge product $\wedge$ of differential forms, i.e.,
for $\omega,\omega_{1}\in\Lambda(N)$, instead of writing $\omega\wedge
\omega_{1}$, we will simply write $\omega\omega_{1}$.

Let $\alpha:A\longrightarrow N$ be an affine bundle (for instance, a vector
bundle) and $F:N_{1}\longrightarrow N$ a smooth map of manifolds. Let
$\mathscr{A}$ be the affine space of smooth sections of $\alpha$. The affine
bundle on $N_{1}$ induced by $\alpha$ via $F$ will be denoted by $F^{\circ
}(\alpha):F^{\circ}(A)\longrightarrow N$:
\[%
\begin{array}
[c]{c}%
\xymatrix{F^\circ(A) \ar[r] \ar[d]_-{F^\circ(\alpha)} & A \ar[d]^-{\alpha} \\ N_1 \ar[r]^-F & N }
\end{array}
,
\]
and the space of its section by $F^{\circ}(\mathscr{A})$. For any section
$a\in\mathscr{A}$ there exists a unique section, which we denote by $F^{\circ
}(a)\in F^{\circ}(\mathscr{A})$, such that the diagram
\[
\xymatrix{F^\circ(A) \ar[r]  & A  \\
N_1 \ar[r]^-F    \ar[u]^-{F^\circ(a)}                    &  N \ar[u]_-{a}}
\]
commutes. If $F:N_{1}\longrightarrow N$ is the embedding of a submanifold, we
also write $\bullet\ |_{F}$ for $F^{\circ}({}\bullet{})$.

We will often understand the sum over repeated upper-lower (multi)indexes. Our
notations about multiindexes are the following. We will use the capital
letters $I,J,K$ for multiindexes. Let $n$ be a positive integer. A multiindex
of length $k$ is a $k$tuple of indexes $I=(i_{1},\ldots,i_{k})$, $i_{1}%
,\ldots,i_{k}\leq n$. We identify multiindexes differing only by the order of
the entries. If $I$ is a multiindex of length $k$, we put $|I|:=k$. Let
$I=(i_{1},\ldots,i_{k})$ and $J=(j_{1},\ldots,j_{l})$ be multiindexes, and $i$
an index. We denote by $IJ$ (resp.{} $Ii$) the multiindex $(i_{1},\ldots
,i_{k},j_{1},\ldots,j_{l})$ (resp.{} $(i_{1},\ldots,i_{k},i)$).

Let $\xi:P\longrightarrow M$ be a fiber bundle. For $0\leq l\leq k\leq\infty$,
we denote by $\xi_{k}:J^{k}\xi\longrightarrow M$ the bundle of $k$-jets of
local sections of $\xi$, and by $\xi_{k,l}:J^{k}\xi\longrightarrow J^{l}\xi$
the canonical projection. For all $k\geq0$, $\xi_{k+1,k}:J^{k+1}%
\xi\longrightarrow J^{k}\xi$ is an affine bundle. For any (local) section
$s:M\longrightarrow P$ of $\pi$, we denote by $j_{k}s:M\longrightarrow
J^{k}\xi$ its $k$th jet prolongation. Let $\ldots,x^{i},\ldots$ be coordinates on $M$ and $\ldots,x^{i},\ldots,y^{a},\ldots$ bundle coordinates on $E$. We denote by $\ldots,x^{i}%
,\ldots,y^{a}{}_{|I},\ldots$ (or simply $\ldots,x^{i},\ldots,y_{I}^{a},\ldots$
if this does not lead to confusion) the associated jet coordinates on
$J^{k}\xi$, $|I|{}\leq k$. Let $\xi^{\prime}:P^{\prime}\longrightarrow M$ be
another fiber bundle and $\Phi:P\longrightarrow P^{\prime}$ a bundle morphism.
For any (local) section $s:M\longrightarrow P$ of $\xi$, $\Phi\circ s:M\longrightarrow P^{\prime}$ is a (local) section of
$\xi^{\prime}$ and, for all $k\geq0$, there exists a unique morphism of bundles over $M$, $\Phi_{\lbrack
k]}:J^{k}\xi\longrightarrow J^{k}\xi^{\prime}$, such that
$\Phi_{\lbrack k]}\circ j_{k}s=j_{k}(\Phi\circ s)$ for any local section $s$
of $\xi$.

Let $\xi$ be as above. We denote by $\mathscr{M}\xi$ and $J^{\dag}\xi$ the
\emph{multimomentum bundle} and the \emph{reduced multimomentum bundle} of $\xi$,
respectively, and by $\tau^{\dag}\xi:J^{\dag}\xi\longrightarrow M$ and
$\tau_{0}^{\dag}\xi:J^{\dag}\xi\longrightarrow P$ the canonical projections (see, for instance, \cite{r05}).
Sections of $\mathscr{M}\xi\longrightarrow P$ are affine maps from sections of
$\xi_{1,0}:J^{1}\xi\longrightarrow P$ to top forms on $M$. Similarly, sections
of $\tau_0^{\dag}\xi$ are the linear parts of sections of $\mathscr{M}\xi
\longrightarrow P$, so that the projection $\mathscr{M}\xi\longrightarrow
J^{\dag}\xi$ is a one dimensional affine bundle.

\section{The \textquotedblleft Technology\textquotedblright\ of
Connections\label{SecConn}}

In this section we discuss some standard and not so standard definitions and
results about (Ehresmann) connections.

Let $\xi:P\longrightarrow M$ and $\ldots,x^{i},\ldots,y^{a},\ldots$ be as
in the previous section. We denote by $\Lambda
_{1}(P,\xi)=\bigoplus_{k}\Lambda_{1}^{k}(P,\xi)\subset\Lambda(P)$ the
differential (graded) ideal in $\Lambda(P)$ made of differential forms on $P$
vanishing when pulled-back to fibers of $\xi$, by $\Lambda_{q}(P,\xi
)=\bigoplus_{k}\Lambda_{q}^{k}(P,\xi)$ its $q$-th exterior power, $q\geq0$,
and by $V\!\Lambda(P,\xi)=\bigoplus_{k}V\!\Lambda^{k}(P,\xi)$ the quotient
differential algebra $\Lambda(P)/\Lambda_{1}(P,\xi)$, $d^{V}:V\!\Lambda
(P,\xi)\longrightarrow V\!\Lambda(P,\xi)$ being its (quotient) differential.
By abusing the notation, we also denote by $d^{V}$ the (quotient) differential
in $\Lambda_{q}(P,\xi)/\Lambda_{q+1}(P,\xi)\simeq V\!\Lambda(P,\xi
)\otimes\Lambda_{q}^{q}(P,\xi)$. There are canonical isomorphisms $\Lambda
_{q}^{q}(P,\xi)\simeq\xi^{\circ}(\Lambda^{q}(M))$, $q\geq0$, and $V\!\Lambda
^{1}(P,\xi)\simeq\operatorname{Hom}(V\!D(P,\xi),C^{\infty}(P))$, $V\!D(P,\xi)$
being the module of $\xi$-vertical vector fields over $P$. In the
following we will understand all the isomorphisms above. 

$\Lambda_{n-1}^{n}(P,\xi)$ and
$V\!\Lambda^{1}(P,\xi)\otimes\Lambda_{n-1}^{n-1}(P,\xi)$ (in this section,
tensor products are over $C^{\infty}(P)$, unless otherwise indicated) identify
canonically with the modules of sections of $\mathscr{M}\xi\longrightarrow P$
and $\tau_{0}^{\dag}\xi:J^{\dag}\xi\longrightarrow P$ respectively.
Accordingly, there is a tautological $n$-form $\Theta$ on $\mathscr{M}\xi$
with the following universal property. For any $\eta\in\Lambda_{n-1}^{n}%
(P,\xi)$, $\eta=\eta^{\ast}(\Theta)$. Similarly, there is a tautological
element $\underline{\Theta}\in V\!\Lambda^{1}(J^{\dag}\xi,\tau^{\dag}%
\xi)\otimes\Lambda_{n-1}^{n-1}(P,\xi)$ with the following universal property.
For any $\underline{\eta}\in$ $V\Lambda^{1}(P,\xi)\otimes\Lambda_{n-1}%
^{n-1}(P,\xi)$, $\underline{\eta}=\underline{\eta}^{\ast}(\underline{\Theta})$.

We denote by $C(P,\xi)$ the affine space of (Ehresmann) connections in
$\xi$. $C(P,\xi)$ identifies canonically with the (affine) space of sections of the
first jet bundle $\xi_{1,0}:J^{1}\xi\longrightarrow P$ and in the following we
will understand such identification. In particular, for $\nabla\in C(P,\xi)$,
we put $\ldots,\nabla_{i}^{a}:=\nabla^{\ast}(y_{i}^{a}),\ldots$, where $\ldots
,y_{i}^{a},\ldots$ are jet coordinates in $J^{1}\xi$. The $\nabla_{i}^{a}$'s
are the \emph{symbols} of the connection $\nabla$. Recall that a (local)
section $\sigma:M\longrightarrow P$ is $\nabla$-constant for some connection
$\nabla\in C(P,\xi)$ iff, by definition, $\nabla\circ\sigma=j_{1}\sigma$, where
$j_{1}\sigma:M\longrightarrow J^{1}\xi$ is the first jet prolongation of
$\sigma$. A connection $\nabla$ in $P$ determines a splittings of the exact
sequence
\begin{equation}
0\longrightarrow V\!D(P,\xi)\longrightarrow D(P)\longrightarrow\xi^{\circ
}(D(M))\longrightarrow0 \label{ExSeq}%
\end{equation}
and its dual%
\begin{equation}
0\longleftarrow V\!\Lambda^{1}(P,\xi)\longleftarrow\Lambda^{1}(P)\longleftarrow
\Lambda_{1}^{1}(P,\xi)\longleftarrow0. \label{ExSeqDual}%
\end{equation}
Thus, using $\nabla$ one can lift a vector field $X$ on $M$ to a vector field
$X^{\nabla}$ transversal to fibers of $\xi$. Moreover, $\nabla$ determines an
isomorphism
\[
\Lambda(P)\simeq\bigoplus_{p,q}V\!\Lambda^{p}(P,\xi)\otimes\Lambda_{q}%
^{q}(P,\xi),
\]
and, in particular, for any $p,q$, a projection
\[
i^{p,q}(\nabla):\Lambda^{p+q}(P)\longrightarrow V\!\Lambda^{p}(P,\xi
)\otimes\Lambda_{q}^{q}(P,\xi),
\]
and an embedding
\[
e^{p,q}(\nabla):V\!\Lambda^{p}(P,\xi)\otimes\Lambda_{q}^{q}(P,\xi
)\longrightarrow\Lambda^{p+q}(P)
\]
taking its values in $\Lambda_{q}^{p+q}(P,\xi)$. For instance, $e^{1,n-1}%
(\nabla)$ is geometrically described by a section $\Sigma_{\nabla}:J^{\dag}%
\xi\longrightarrow\mathscr{M}\xi$ of the affine bundle $\mathscr{M}\xi
\longrightarrow J^{\dag}\xi$.

Finally, every connection $\nabla$ defines a vector valued differential
$2$-form, its curvature, $R^{\nabla}\in\Lambda_{2}^{2}(P,\xi)\otimes
V\!D(P,\xi)$, via
\[
R^{\nabla}(X,Y):=[X^{\nabla},Y^{\nabla}]-[X,Y]^{\nabla}, \quad X,Y\in D(M)
\]
Locally,
\[
R^{\nabla}=R_{ij}^{a}dx^{i}dx^{j}\otimes\tfrac{\partial}{\partial y^{a}},\quad
R_{ij}^{a}=\tfrac{1}{2}(D_{i}\nabla_{j}^{a} - D_{j}\nabla_{i}^{a})\circ\nabla.
\]
where $D_{i}:=\partial_{i}+y_{i}^{a}\tfrac{\partial}{\partial y^{a}}$,
$\partial_{i}:=\tfrac{\partial}{\partial x^{i}}$, $i=1,\ldots,n$. A connection
$\nabla$ is flat iff, by definition, $R^{\nabla}=0$. If $\nabla$ is a flat
connection in $\xi$, then $P$ is locally foliated by (local) $\nabla$-constant
sections of $\xi$.

In the following, it will be also useful the concept of a \emph{relative
connection}. Let $\xi:P\longrightarrow M$ be as above, $\zeta:N\longrightarrow
M$ another fiber bundle, and $F:N\longrightarrow P$ a bundle morphism. A
\emph{relative connection along }$F$ is an element $\square$ of the affine
space $F^{\circ}(C(P,\xi))$, i.e., a section of the induced bundle $F^{\circ
}(\xi_{1,0}):F^{\circ}(J^{1}\xi)\longrightarrow N$, or, which is the same, a
map $\square:N\longrightarrow J^{1}\xi$ such that $\xi_{1,0}\circ\square=F$. A
(local) section $\Sigma:M\longrightarrow N$ of $\zeta$ is $\square$-constant,
for a relative connection $\square\in F^{\circ}(C(P,\xi))$, iff, by
definition, $\square\circ\Sigma=j_{1}(F\circ\Sigma)$. A relative connection
$\square$ along $F$ determines a splittings of the exact sequence
\[
0\longrightarrow F^{\circ}(V\!D(P,\xi))\longrightarrow F^{\circ}%
(D(P))\longrightarrow\zeta^{\circ}(D(M))\longrightarrow0
\]
and its dual%
\[
0\longleftarrow F^{\circ}(V\!\Lambda^{1}(P,\xi))\longleftarrow F^{\circ}%
(\Lambda^{1}(P))\longleftarrow\Lambda_{1}^{1}(N,\zeta)\longleftarrow0,
\]
which are obtained from sequences (\ref{ExSeq}) and (\ref{ExSeqDual}) by
tensorizing for $C^{\infty}(N)$. Thus, using $\square$ one can lift a vector
field $X$ on $M$ to a relative vector field $X^{\square}$ along $F$ transversal
to fibers of $\xi$. Moreover, $\square$ determines an isomorphism
\[
F^{\circ}(\Lambda(P))\simeq\bigoplus_{p,q}F^{\circ}(V\!\Lambda^{p}%
(P,\xi))\otimes\Lambda_{q}^{q}(N,\zeta),
\]
and, in particular, for any $p,q$, a projection
\[
i^{p,q}(\square):F^{\circ}(\Lambda^{p+q}(P))\longrightarrow F^{\circ
}(V\!\Lambda^{p}(P,\xi))\otimes\Lambda_{q}^{q}(N,\zeta),
\]
and an embedding
\[
e^{p,q}(\square):F^{\circ}(V\!\Lambda^{p}(P,\xi))\otimes\Lambda_{q}%
^{q}(N,\zeta)\longrightarrow F^{\circ}(\Lambda^{p+q}(P))
\]
taking its values in $F^{\circ}(\Lambda_{q}^{p+q}(P,\xi))$. For instance,
$e^{1,n-1}(\square)$ is geometrically described by a section $\Sigma_{\nabla
}:F^{\circ}(J^{\dag}\xi)\longrightarrow F^{\circ}(\mathscr{M}\xi)$ of the
affine bundle $F^{\circ}(\mathscr{M}\xi)\longrightarrow F^{\circ}(J^{\dag}%
\xi)$.

\begin{example}
Let $\xi:P\longrightarrow M$ be as above, $y\in P$ and $z\in\xi_{1,0}%
^{-1}(y)\subset J^{1}\xi$. $z$ can be understood as a relative connection
along the embedding $y:\ast\longrightarrow P$ at $y$ of the one point manifold
$\ast$. If $\omega\in\Lambda_{n-1}^{n+1}(P,\xi)$ is a \emph{PD-Hamiltonian
system} on $\xi$ in the sense of \cite{v09b}, then its \emph{first constraint
submanifold }$\mathscr{P}$ is defined as
\[
\mathscr{P}:=\{y\in P\;|\;i^{1,n}(z)(\omega|_{y})=0\text{ for some }z\in
\xi_{1,0}^{-1}(y)\}\subset P.
\]
\end{example}

\begin{example}
Let $\xi:P\longrightarrow M$ be as above and $\sigma:M\longrightarrow P$ a
(local) section of $\xi$. It is sometimes useful to understand $j_{1}%
\sigma:M\longrightarrow J^{1}\xi$ as a relative connection along $\sigma$. For
instance, if $\omega\in\Lambda_{n-1}^{n+1}(P,\xi)$ is a PD-Hamiltonian system
on $\xi$, the \emph{PD-Hamilton equations} \cite{v09b} for $\sigma$ read
\[
i^{1,n}(j_{1}\sigma)(\omega|_{\sigma})=0.
\]
\end{example}

\begin{example}
\label{Examp}Consider a fiber bundle $\pi:E\longrightarrow M$ and, for some
$l$, the projection $\pi_{l+1,l}:J^{l+1}\pi\longrightarrow J^{l}\pi$. Let $\ldots,x^i,\ldots$ be coordinates on $M$ and $\ldots,x^i,\ldots,u^\alpha,\ldots$ bundle coordinates on $E$. There is
a canonical relative connection along $\pi_{l+1,l}$. Namely, recall that
$J^{l+1}\pi$ is canonically embedded into $J^{1}\pi_{l}$ via
\[
\mathscr{C}:J^{l+1}\pi\ni(j_{l+1}s)(x)\longmapsto(j_{1}(j_{l}s))(x)\in
J^{1}\pi_{l},
\]
where $s$ is a local section of $\pi$ and $x\in M$. $\mathscr{C}$ is a
relative connection and locally $\ldots,\mathscr{C}^{\ast}(u_{I}^{\alpha}{}%
_{|i})=u_{Ii}^{\alpha},\ldots$, $|I|{}\leq l$. $\mathscr{C}$ will be called the \emph{Cartan
(relative) connection along }$\pi_{l+1,l}$. Notice that $\mathscr{C}$-constant
sections of $\pi_{l+1}$ are precisely holonomic sections, i.e., sections of
the form $j_{l+1}s:M\longrightarrow J^{l+1}\pi$ for some section $s$ of $\pi$.
Finally one can draw the following commutative diagram
\begin{equation}%
\begin{array}
[c]{c}%
%TCIMACRO{\TeXButton{TeX field}{\xymatrix@C=2pt@R=15pt{ & \pi_{l+1,l}%
%^\circ(J^\dag\pi_l) \ar[ddrrr]^-{h} \ar[ddl] \ar[ddr]   \ar@
%/^1.2pc/[rr]^-{\Sigma_\mathscr{C}} & & \pi_{l+1,l}^\circ(\mathscr{M}\pi
%_l) \ar[ddr] \ar[ll] \ar[ddlll]|!{[ll];[ddr]}\hole|!{[ll];[ddl]}\hole& \\
%& & & &                                                                                                             \\
%J^{l+1} \ar[ddr] & & J^\dag\pi_l \ar[ddl]   & & \mathscr{M}\pi_l \ar
%[ll] \ar
%[ddlll]                                                                   \\
%& & & &
%\\ & J^l \ar[dd]^-{\pi_l} &  & & \\ & & & & \\ & M & & &}}}%
%BeginExpansion
\xymatrix@C=2pt@R=15pt{ & \pi_{l+1,l}^\circ(J^\dag\pi_l) \ar[ddrrr]^-{h}
\ar[ddl] \ar[ddr]   \ar@/^1.2pc/[rr]^-{\Sigma_\mathscr{C}} & & \pi
_{l+1,l}^\circ(\mathscr{M}\pi_l) \ar[ddr] \ar[ll] \ar[ddlll]|!{[ll];[ddr]}%
\hole|!{[ll];[ddl]}\hole& \\
& & & &                                                                                                             \\
J^{l+1} \ar[ddr] & & J^\dag\pi_l \ar[ddl]   & & \mathscr{M}\pi_l \ar
[ll] \ar
[ddlll]                                                                   \\
& & & &
\\ & J^l \ar[dd]^-{\pi_l} &  & & \\ & & & & \\ & M & & &}%
%EndExpansion
\end{array}
\label{Diag}%
\end{equation}
where $h$ is the composition of $\Sigma_{\mathscr{C}}$ and the projection
$\pi_{l+1,l}^{\circ}(\mathscr{M}\pi_{l})\longrightarrow\mathscr{M}\pi_{l}$.
Notice that $h$ is a morphism of bundles over $J^{\dag}\pi_{l}$. Moreover, let
$\ldots,p_{\alpha}^{I.i},\ldots,p$ be standard coordinates on
$\mathscr{M}\pi_{l}$ \cite{r05} associated to jet coordinates $\ldots,x^{i},\ldots,u_{I}^{\alpha},\ldots
$ on $J^l$, $|I|{}\leq l$. Then $h$ is locally given by $h^{\ast}(p)=\sum
_{|I|{}\leq l}p_{\alpha}^{I.i}u_{Ii}^{\alpha}$.
\end{example}

\section{Lagrangian-Hamiltonian Formalism\label{SecLHForm}}

In this section we present a finite dimensional version of the formalism of
\cite{v09}.

Let $\pi:E\longrightarrow M$ be a fiber bundle and $\ldots,x^{i},\ldots
,u_{I}^{\alpha},\ldots,p_{\alpha}^{I.i},\ldots,p$ coordinates on $\mathscr{M}\pi_l$ as in Example
\ref{Examp}, $i=1,\ldots,n$, $\alpha=1,\ldots,m$, $\dim M=n$, $\dim E=m+n$. In
the following we put $J^{k}:=J^{k}\pi$, $k\geq0$.

\begin{definition}
A \emph{Lagrangian field theory} of order $l+1$ is a pair $(\pi,\mathscr{L})$
where $\mathscr{L}$ is a \emph{Lagrangian density} of order $l+1$,
i.e., an element of $\Lambda_{n}^{n}(J^{l+1},\pi_{l+1})$.
\end{definition}

Let $(\pi,\mathscr{L})$ be a Lagrangian field theory of order $l+1$.
$\mathscr{L}$ is locally given by $\mathscr{L}=Ld^{n}x$, $d^n x := dx^1 \cdots dx^n$, $L$ being a local
function on $J^{l+1}$. Extremals of the variational principle $\int
\mathscr{L}$ are solutions of the \emph{Euler-Lagrange equations}, i.e.,
(local) sections of $\pi$, $s:M\longrightarrow E$ such that
\begin{equation}
(j_{2l+2}s)^{\circ}(\boldsymbol{E}(\mathscr{L}))=0, \label{ELE}%
\end{equation}
$\boldsymbol{E}(\mathscr{L})\in\pi_{2l+2,0}^{\circ}(V\!\Lambda(E,\pi
)\otimes_{A}\Lambda_{n}^{n}(M))$ being the Euler-Lagrange form (see, for
instance, \cite{b...99}), which is locally given by
\[
\boldsymbol{E}(\mathscr{L})=\tfrac{\delta L}{\delta u^{\alpha}}d^{V}%
\!u^{\alpha}\otimes d^{n}x
\]
where $\tfrac{\delta L}{\delta u^{\alpha}}:=\sum_{|I|{}\leq l+1}(-)^{|I|}%
D_{I}\partial_{\alpha}^{I}L$, $\alpha=1,\ldots,m$, are the variational
derivatives of $L$, $D_{(i_{1},\ldots,i_{k})}:=D_{i_{1}}\circ\cdots\circ
D_{i_{k}}$, $D_{i}=\partial_{i}+\sum_{|I|{}\geq0}u_{Ii}^{\alpha}%
\partial_{\alpha}^{I}$ is the $i$th total derivative, $i=1,\ldots,n$, and
$\partial_{\alpha}^{I}:=\tfrac{\partial}{\partial u_{I}^{\alpha}}$,
$\alpha=1,\ldots,m$, $|I|{}\geq0$.

Consider the commutative subdiagram of (\ref{Diag})
\[%
\begin{array}
[c]{c}%
%TCIMACRO{\TeXButton{TeX field}{\xymatrix{ \pi_{l+1,l}^\circ(J^\dag
%)  \ar[r]^-{h} \ar[d]_-{\mathfrak{p}} & \mathscr{M}\ar[d] \\
%J^{l+1} \ar[r]^-{\pi_{l+1,l}} & J^l
%}}}%
%BeginExpansion
\xymatrix{ \pi_{l+1,l}^\circ(J^\dag)  \ar[r]^-{h} \ar[d]_-{\mathfrak{p}}
& \mathscr{M}\ar[d] \\
J^{l+1} \ar[r]^-{\pi_{l+1,l}} & J^l
}%
%EndExpansion
\end{array}
.
\]
where $\mathfrak{p}$ is the canonical projection. Here and in what follows
$J^{\dag}:=J^{\dag}\pi_{l}$ and $\mathscr{M}:=\mathscr{M}\pi_{l}$. Similarly,
we put $\tau_{0}^{\dag}:=\tau_{0}^{\dag}\pi_{l}:J^{\dag}\longrightarrow J^{l}$
and $\tau^{\dag}:=\tau^{\dag}\pi_{l}:J^{\dag}\longrightarrow M$. Consider also
the $n$-form $\theta^{\prime}$ on $\pi_{l+1,l}^{\circ}(J^{\dag})$
defined as
\[
\theta^{\prime}:=h^{\ast}(\Theta)+\mathfrak{p}^{\ast}(\mathscr{L}),
\]
$\Theta$ being the tautological $n$-form on $\mathscr{M}$. $\theta^{\prime}$
is locally given by
\[
\theta^{\prime}=%
%TCIMACRO{\tsum \nolimits_{|I|{}\leq l}}%
%BeginExpansion
{\textstyle\sum\nolimits_{|I|{}\leq l}}
%EndExpansion
p_{\alpha}^{I.i}du_{I}^{\alpha}d^{n-1}x_{i}-Ed^{n}x,
\]
where $d^{n-1}x_i := i_{\partial_i} d^n x$ and $E:=\sum_{|I|{}\leq l}p_{\alpha}^{I.i}u_{Ii}^{\alpha}-L$. $d\theta
^{\prime}$ is a PD-Hamiltonian system on the bundle $\pi_{l+1,l}^{\circ
}(J^{\dag})\longrightarrow M$, whose first constraint submanifold
$\mathscr{P}$ is locally defined by
\begin{equation}
\partial_{\alpha}^{I}L-%
%TCIMACRO{\tsum \nolimits_{|J|{}\leq l}}%
%BeginExpansion
{\textstyle\sum\nolimits_{|J|{}\leq l}}
%EndExpansion
\delta_{Ji}^{I}p_{\alpha}^{J.i}=0,\quad|I|{}=l+1, \label{EqP}%
\end{equation}
where the symbol $\delta^I_K$ is equal to $1$ if the mutliindexes $I$ and $K$ coincide and is equal to $0$ otherwise. In particular, $\mathscr{P}$ has the same dimension as $J^{\dag}$ and
$\mathscr{P}\longrightarrow J^{l+1}$ is an affine bundle. However the map
$\mathscr{P}\longrightarrow J^{\dag}$ needs not be submersive. As will be
clear in a moment, it is natural to give the following

\begin{definition}
The theory $(\pi,\mathscr{L})$ is \emph{hyperregular }iff the map
$\mathscr{P}\longrightarrow J^{\dag}$ is a diffeomorphism.
\end{definition}

The above definition generalizes the standard definition of a hyperregular
Lagrangian system on the tangent bundle of a configuration manifold. In the
following we will always assume $(\pi,\mathscr{L})$ to be hyperregular. In
particular, the matrix%
\begin{equation}
\mathbf{H}:=\left\Vert (\partial_{\beta}^{K}\partial_{\alpha}^{I}%
L)(\theta)\right\Vert {}_{(\beta,K)}^{(\alpha,I)},\quad|I|,|K|{}=l+1,
\label{Hessian}%
\end{equation}
where the pairs $(\alpha,I)$ and $(\beta,K)$ are understood as single indexes,
has maximum rank at every point $\theta\in J^{l+1}$. The case when matrix
(\ref{Hessian}) has lower rank (which is physically relevant for gauge
theories) will be treated elsewhere.

Notice that, inverting the diffeomorphism $\mathscr{P}\longrightarrow J^{\dag
}$ (and composing with $i_{\mathscr{P}}:\mathscr{P}\longrightarrow\pi
_{l+1,l}^{\circ}(J^{\dag})$) we get a section $\mathfrak{s}:J^{\dag
}\longrightarrow\pi_{l+1,l}^{\circ}(J^{\dag})$ of $\pi_{l+1,l}^{\circ}%
(J^{\dag})\longrightarrow J^{\dag}$. Put $\theta:=\mathfrak{s}^{\ast}%
(\theta^{\prime})$. Then, locally
\[
\theta=%
%TCIMACRO{\tsum \nolimits_{|I|{}\leq l}}%
%BeginExpansion
{\textstyle\sum\nolimits_{|I|{}\leq l}}
%EndExpansion
p_{\alpha}^{I.i}du_{I}^{\alpha}d^{n-1}x_{i}-Hd^{n}x,
\]
where $H:=\mathfrak{s}^{\ast}(E)$. Moreover, $\omega:=d\theta$ is a
PD-Hamiltonian system on $J^{\dag}$ and determines PD-Hamilton equations
\begin{equation}
i^{1,n}(j_{1}\sigma)(\omega|_{\sigma})=0 \label{HDWE}%
\end{equation}
for (local) sections $\sigma$ of $\pi_{l}$. Equations (\ref{HDWE}) read
locally
\begin{equation}
\left\{
\begin{array}
[c]{l}%
p_{\alpha}^{I.i},_{i}=-\tfrac{\partial H}{\partial u_{I}^{\alpha}}\\
u_{I}^{\alpha},_{i}=\tfrac{\partial H}{\partial p_{\alpha}^{I.i}}%
\end{array}
\right.  ,\quad|I|{}\leq l, \label{DEloc}%
\end{equation}
where a \textquotedblleft$\bullet,_{i}$\textquotedblright\ denotes
differentiation of \textquotedblleft$\bullet$\textquotedblright\ with respect
to $x^{i}$, $i=1,\ldots,n$. Eqs.{} (\ref{DEloc}) are higher order de Donder field
equations \cite{d35}. Accordingly, Equations (\ref{HDWE}) will be refereed to
as the \emph{Hamilton-de Donder-Weyl Equations} \emph{(HDWE) }determined by
the field theory $(\pi,\mathscr{L})$.

In the remaining part of this section we provide an alternative description of
$\theta$. First of all, notice that, in view of the universal property of
$\Theta$, there exists a unique morphism of bundles over $J^{l}$,
$\mathbb{L}:\pi_{l+1,l}^{\circ}(J^{\dag})\longrightarrow\mathscr{M}$, such
that%
\[
\mathfrak{p}^{\ast}(\mathscr{L})=\mathbb{L}^{\ast}(\Theta).
\]
$\mathbb{L}$ is locally given by $\ldots,\mathbb{L}^{\ast}(p_{\alpha}%
^{I.i})=0,\ldots,\mathbb{L}^{\ast}(p)=-L$, $|I|{}\leq l$.

Consider the section $\mathscr{H}:=\mathfrak{s}\circ(h+\mathbb{L}):J^{\dag
}\longrightarrow\mathscr{M}$ of $\mathscr{M}\longrightarrow J^{\dag}$.

\begin{definition}
$\mathscr{H}$ is called the \emph{Hamiltonian section} of the Lagrangian field theory
$(\pi,\mathscr{L})$.
\end{definition}

\begin{proposition}
$\theta=\mathscr{H}^{\ast}(\Theta)$.
\end{proposition}

\begin{proof}
Indeed,
\begin{align*}
\theta &  =\mathfrak{s}^{\ast}(\theta^{\prime})\\
&  =\mathfrak{s}^{\ast}(h^{\ast}(\Theta)+\mathfrak{p}^{\ast}(\mathscr{L}))\\
&  =\mathfrak{s}^{\ast}(h^{\ast}(\Theta)+\mathbb{L}^{\ast}(\Theta))\\
&  =(\mathfrak{s}^{\ast}\circ(h+\mathbb{L})^{\ast})(\Theta)\\
&  =((h+\mathbb{L})\circ\mathfrak{s})^{\ast}(\Theta)\\
&  =\mathscr{H}^{\ast}(\Theta).
\end{align*}
\end{proof}

It follows from the above proposition that locally $\mathscr{H}^{\ast}(p)=H$.

\section{The Inverse Legendre Transform\label{SecInvLeg}}

Notice that, in the previous section, we introduced the Hamiltonian formalism
for the Lagrangian theory $(\pi,\mathscr{L})$ without defining a (higher order
analogue of the) Legendre transform. Actually, a Legendre transform can be
only introduced by means of an external structure such as a coordinate system,
a Lepagean equivalent of $\mathscr{L}$ \cite{k84}, a connection in $\pi$
\cite{k84b}, or a Legendre form \cite{sc90,av04}. Among these structures,
there is generically no distinguished one. Therefore, for the sake of the
naturality of the formalism, we prefer not to introduce any Legendre
transform. Nevertheless, the concept of an \textquotedblleft\emph{inverse
Legendre transform}\textquotedblright\ can be introduced without ambiguity as follows.

Put $\mathsf{leg}^{-1}:=\mathfrak{p}\circ\mathfrak{s}:J^{\dag}\longrightarrow
J^{l+1}$. Clearly, $\mathsf{leg}^{-1}$ is a(n affine) bundle and a morphism of
bundles over $J^{l}$.

\begin{definition}
$\mathsf{leg}^{-1}:J^{\dag}\longrightarrow J^{l+1}$ is called the
\emph{inverse Legendre transform}.
\end{definition}

In the remaining part of this section we provide an alternative description of
the inverse Legendre transform.

First of all, notice that $\omega$ is an unconstrained PD-Hamiltonian system
on $J^{\dag}$ in the sense that its first constraint submanifold coincides with
the whole $J^{\dag}$. Namely, for any $P\in J^{\dag}$ the subset
\[
\operatorname{Ker}\omega|_{P}:=\{z\in(\tau^{\dag})_{1,0}^{-1}(P)\;|\;i^{1,n}%
(z)(\omega|_{P})=0\}\subset(\tau^{\dag})_{1,0}^{-1}(P)
\]
is non-empty.

\begin{proposition}
\label{Prop2}For any $P\in J^{\dag}$, $(\tau_{0}^{\dag})_{[1]}(z)\in J^{1}%
\pi_{l}$ is independent of the choice of $z\in$ $\operatorname{Ker}\omega
|_{P}$.
\end{proposition}

\begin{proof}
Let $P\in J^{\dag}$ and $z\in\operatorname{Ker}\omega|_{P}$. Choose standard
coordinates $\ldots,u_{I}^{\alpha},\dots,p_{\alpha}^{I.i},\ldots$ around $P$
and let $\ldots,z_{I}^{\alpha}{}_{|j}:=u_{I}^{\alpha}{}_{|j}(z),\ldots
,z_{\alpha}^{I.i}{}_{|j}:=p_{\alpha}^{I.i}{}_{|j}(z),\ldots$. Then
\[
0=i^{1,n}(z)(\omega|_{P})=-[%
%TCIMACRO{\tsum \nolimits_{|I|{}\leq l}}%
%BeginExpansion
{\textstyle\sum\nolimits_{|I|{}\leq l}}
%EndExpansion
(z_{\alpha}^{I.i}{}_{|i}+\tfrac{\partial H}{\partial u_{I}^{\alpha}}%
(P))d^{V}\!u_{I}^{\alpha}+%
%TCIMACRO{\tsum \nolimits_{|I|{}\leq l}}%
%BeginExpansion
{\textstyle\sum\nolimits_{|I|{}\leq l}}
%EndExpansion
(z_{I}^{\alpha}{}_{|i}-\tfrac{\partial H}{\partial p_{\alpha}^{I.i}}%
(P))d^{V}\!p_{\alpha}^{I.i}]\otimes d^{n}x|_{P}.
\]
In particular,
\[
z_{I}^{\alpha}{}_{|i}=\tfrac{\partial H}{\partial p_{\alpha}^{I.i}}%
(P),\quad|I|{}\leq l,
\]
which do only depend on $P$. Since the $z_{I}^{\alpha}{}_{|i}$'s completely
identify $(\tau_{0}^{\dag})_{[1]}(z)$, the assertion follows.
\end{proof}

In view of the above proposition (and its proof), the map
\[
\mathbb{F\!}\mathscr{H}:J^{\dag}\ni P\longmapsto(\tau_{0}^{\dag})_{[1]}(z)\in
J^{1}\pi_{l},\quad z\in\operatorname{Ker}\omega|_{P},
\]
is a well-defined morphism of bundles over $J^{l}$ locally given by
$\ldots,\mathbb{F\!}\mathscr{H}^{\ast}(u_{I}^{\alpha}{}_{|i})=\tfrac{\partial
H}{\partial p_{\alpha}^{I.i}},\ldots$, $|I|{}\leq l$. In particular, $\mathbb{F\!}%
\mathscr{H}$ identifies with the \emph{fiber derivative} of $\mathscr{H}$ (see
for instance \cite{fr05}).

\begin{proposition}
Diagram
\[%
\begin{array}
[c]{c}%
%TCIMACRO{\TeXButton{TeX field}{\xymatrix{ J^1 \pi_l & & J^\dag\pi
%_l \ar[ll]_-{\mathbb{F}\!\mathscr{H}} \ar[lld]^-{\mathsf{leg}^{-1}}  \\
%J^{l+1} \ar[u]^-{\mathscr{C}}  &  &
%}}}%
%BeginExpansion
\xymatrix{ J^1 \pi_l & & J^\dag\pi_l \ar[ll]_-{\mathbb{F}\!\mathscr{H}}
\ar[lld]^-{\mathsf{leg}^{-1}}  \\
J^{l+1} \ar[u]^-{\mathscr{C}}  &  &
}%
%EndExpansion
\end{array}
\]
commutes.
\end{proposition}

\begin{proof}
Put $\ldots,\mathfrak{s}_{J}^{\alpha}:=\mathfrak{s}^{\ast}(u_{J}^{\alpha})\ldots$,
$|J|{}\leq l+1$. Clearly, $\ldots,\mathfrak{s}_{J}^{\alpha}=u_{J}^{\alpha},\ldots$ for
$|J|{}<l+1$, while, for $|J|{}=l+1$, the $\mathfrak{s}_{J}^{\alpha}$'s are
implicitly defined by
\begin{equation}
(\partial_{\alpha}^{J}L-%
%TCIMACRO{\tsum \nolimits_{|K|{}\leq l}}%
%BeginExpansion
{\textstyle\sum\nolimits_{|K|{}\leq l}}
%EndExpansion
\delta_{Ki}^{J}p_{\alpha}^{K.i})\circ\mathfrak{s}=0,\quad|J|{}=l+1,
\label{ImpEq}%
\end{equation}
which expresses the fact that $\mathfrak{s}$ takes its values in
$\mathscr{P}$. Now, for $|I|{}\leq l$, compute
\begin{align*}
\mathbb{F\!}\mathscr{H}^{\ast}(u_{I}^{\alpha}{}_{|i})  &  =\tfrac{\partial
H}{\partial p_{\alpha}^{I.i}}\\
&  =\tfrac{\partial}{\partial p_{\alpha}^{I.i}}\mathfrak{s}^{\ast}(E)\\
&  =\tfrac{\partial}{\partial p_{\alpha}^{I.i}}(E\circ\mathfrak{s})\\
&  =\tfrac{\partial E}{\partial p_{\alpha}^{I.i}}\circ\mathfrak{s}+%
%TCIMACRO{\tsum \nolimits_{|J|{}\leq l+1}}%
%BeginExpansion
{\textstyle\sum\nolimits_{|J|{}\leq l+1}}
%EndExpansion
(\tfrac{\partial}{\partial p_{\alpha}^{I.i}}\mathfrak{s}_{J}^{\alpha
})\,(\partial_{\alpha}^{J}E\circ\mathfrak{s})\\
&  =u_{Ii}^{\alpha}\circ\mathfrak{s}-%
%TCIMACRO{\tsum \nolimits_{|J|{}=l+1}}%
%BeginExpansion
{\textstyle\sum\nolimits_{|J|{}=l+1}}
%EndExpansion
(\tfrac{\partial}{\partial p_{\alpha}^{I.i}}\mathfrak{s}_{J}^{\alpha
})\,(\partial_{\alpha}^{J}L-{\textstyle\sum\nolimits_{|K|{}\leq l}}\delta_{Ki}^{J}p_{\alpha}^{K.i})\circ
\mathfrak{s}\\
&  =\mathfrak{s}^{\ast}(u_{Ii}^{\alpha})\\
&  =(\mathfrak{s}^{\ast}\circ\mathfrak{p}^{\ast})(u_{Ii}^{\alpha})\\
&  =(\mathfrak{s}^{\ast}\circ\mathfrak{p}^{\ast}\circ\mathscr{C}^{\ast}%
)(u_{I}^{\alpha}{}_{|i})\\
&  =(\mathscr{C}\circ\mathfrak{p}\circ\mathfrak{s})^{\ast}(u_{I}^{\alpha}%
{}_{|i})\\
&  =(\mathscr{C}\circ\mathsf{leg}^{-1})^{\ast}(u_{I}^{\alpha}{}_{|i}).
\end{align*}
Since both $\mathbb{F\!}\mathscr{H}$ and $\mathscr{C}\circ\mathsf{leg}^{-1}$
are morphisms of bundles over $J^{l}$, this concludes the proof.
\end{proof}

In particular, $\operatorname{im}\mathbb{F\!}\mathscr{H}\subset
\operatorname{im}\mathscr{C}$ and $\mathsf{leg}^{-1}$ is obtained from
$\mathbb{F\!}\mathscr{H}$ restricting the codomain to $\operatorname{im}%
\mathscr{C}$.

Finally, notice that Eqs.{} (\ref{HDWE}) cover Euler-Lagrange equations
(\ref{ELE}) via $\mathsf{leg}^{-1}$ in the following sense. If $\sigma
:M\longrightarrow J^{\dag}$ is a solution of (\ref{HDWE}), then $\mathsf{leg}%
^{-1}\circ\sigma=j_{l+1}s:M\longrightarrow J^{l+1}$ for a solution
$s:M\longrightarrow E$ of (\ref{ELE}) (see \cite{v09} for a detailed proof).

\section{Generalized Hamilton-Jacobi Formalism\label{SecGenHJForm}}

In this section we present the analogue of Hamilton-Jacobi formalism for
higher order Lagrangian field theories. In the spirit of \cite{c...06} we
first formulate a generalized Hamilton-Jacobi problem which does only depend
on the field equations (and not directly on the Hamiltonian section).

Let $(\pi,\mathscr{L})$ be a Lagrangian field theory of the order $l+1$. We
use here the same notations as in the previous section. Moreover, we identify
$J^{l+1}$ with the image of the embedding $\mathscr{C}:J^{l+1}\longrightarrow
J^{1}\pi_{l}$, thus understanding $J^{l+1}$ as a submanifold in $J^{1}\pi_{l}$.

\begin{problem}
[generalized Hamilton-Jacobi]\label{HJProb}the generalized
\emph{Hamilton--Jacobi (HJ in the following) problem} for the Lagrangian
theory $(\pi,\mathscr{L})$ consists in finding a section $T$ of $\tau
_{0}^{\dag}:J^{\dag}\longrightarrow J^{l}$ and a flat connection $\nabla$ in
$\pi_{l}:J^{l}\longrightarrow M$ such that%
\begin{equation}
T\circ\gamma\text{ is a solution of HDWE (\ref{HDWE}) for every }%
\nabla\text{-constant section }\gamma:M\longrightarrow J^{l}\text{ of }\pi
_{l}\text{,} \tag{$\star$}\label{HJCond}%
\end{equation}
\end{problem}

Notice that condition (\ref{HJCond}) also implies that every $\nabla$-constant
section is the $l$th jet prolongation of a solution of the Euler-Lagrange
equations. In order to characterize pairs $(T,\nabla)$ satisfying conditions
(\ref{HJCond}) of Problem \ref{HJProb} we put
\[
Z:=%
%TCIMACRO{\tbigcup \nolimits_{P\in J^{\dag}}}%
%BeginExpansion
{\textstyle\bigcup\nolimits_{P\in J^{\dag}}}
%EndExpansion
\operatorname{Ker}\omega|_{P}\subset J^{1}\tau^{\dag}.
\]
$Z$ is the \emph{geometric portrait of the }HDWE in the following sense. A section $\sigma$ of
$\tau^{\dag}$ is a solution of the HDWE iff $\operatorname{im}j_{1}%
\sigma\subset Z$. 

A submanifold $N\subset J^{\dag}$ will be said
$Z$\emph{-compatible} iff for every $P\in N$, there exists $z\in Z$ such that
$(\tau^{\dag})_{1,0}(z)=P$ and $z\subset T_{P}N$.

Let $T$ be a section of $\tau_{0}^{\dag}$. In the following we put $\nabla
^{T}:=\mathsf{leg}^{-1}\circ T$.

\begin{theorem}
\label{Th1}Let $T$ be a section of $\tau_{0}^{\dag}$ and $\nabla$ a flat
connection in $\pi_{l}$. The following assertions are equivalent

\begin{enumerate}
\item $(T,\nabla)$ satisfies condition (\ref{HJCond});

\item $\operatorname{im}(T_{[1]}\circ\nabla)\subset Z$;

\item $\operatorname{im}T$ is $Z$-compatible and $\nabla=\nabla^{T}$.
\end{enumerate}
\end{theorem}

Before proving Theorem \ref{Th1} we prove the following

\begin{lemma}
\label{Lem1}Let $T$ be a section of $\tau_{0}^{\dag}$ and $\nabla$ a (possibly
non flat) connection in $\pi_{l}$. If $\operatorname{im}(T_{[1]}\circ
\nabla)\subset Z$ then $(T,\nabla)$ satisfies condition (\ref{HJCond}) and
$\nabla=\nabla^{T}$.
\end{lemma}

\begin{proof}
Let $T$ and $\nabla$ be as in the hypothesis of the lemma. $T:J^{l}%
\longrightarrow J^{\dag}$ may be understood as a morphism of bundles over $M$.
Let $T_{[1]}:J^{1}\pi_{l}\longrightarrow J^{1}\tau^{\dag}$ be its first jet
prolongation. For every $\nabla$-constant section $\gamma$ we have
\[
\operatorname{im}j_{1}(T\circ\gamma)=\operatorname{im}(T_{[1]}\circ
j_{1}\gamma)\subset\operatorname{im}(T_{[1]}\circ\nabla)\subset Z,
\]
i.e., $T\circ\gamma$ is a solution of HDWE. Moreover, $T_{[1]}$ is
a section of the bundle $(\tau_{0}^{\dag})_{[1]}:J^{1}\tau^{\dag
}\longrightarrow J^{1}\pi_{l}$. Therefore, for any $y\in J^{l}$,
\[
\nabla(y)=((\tau_{0}^{\dag})_{[1]}\circ T_{[1]}\circ\nabla)(y)=\mathbb{F}\!%
\mathscr{H}(T(y))=(\mathsf{leg}^{-1}\circ T)(y).
\]
\end{proof}

Notice that $(T,\nabla^{T})$ determines a section of the bundle $\pi
_{l+1,l}^{\circ}(J^{\dag})\longrightarrow J^{l}$ with values in $\mathscr{P}$
and vice versa.

\begin{proof}
[Proof of Theorem \ref{Th1}]Lemma \ref{Lem1} shows, in particular, that (2)
$\Longrightarrow$ (1).

Now, let $(T,\nabla)$ satisfy condition (\ref{HJCond}), $y\in J^{l}$ and
$x=\pi_{l}(y)$. Since $\nabla$ is flat there is a $\nabla$-constant (local) section
$\gamma$ such that $y=\gamma(x)$. Put $z:=j_{1}(T\circ\gamma)(y)$. Then
$(\tau^{\dag})_{1,0}(z)=(T\circ\gamma)(x)=T(y)$, $z\subset T_{T(y)}%
\operatorname{im}T$ and, since $T\circ\gamma$ is a solution of HDWE, $z\in Z$.
This shows that $\operatorname{im}T$ is $Z$-compatible. Thus (1)
$\Longrightarrow$ (3).

Finally, let $\operatorname{im}T$ be $Z$-compatible. For $y\in J^{l}$, let
$z(y)\in Z$ be such that $(\tau^{\dag})_{1,0}(z(y))=T(y)$ and $z(y)\subset
T_{T(y)}\operatorname{im}T$. Consider the map $\nabla^{\prime}:J^{l}\ni
y\longmapsto(\tau_{0}^{\dag})_{[1]}(z(y))\in J^{1}\pi_{l}$. In view of
Proposition \ref{Prop2}, $\nabla^{\prime}=\nabla^{T}$ and, since $T$ is an
immersion, $(T_{[1]}\circ\nabla^{T})(y)=z(y)\in Z$. Thus (3) $\Longrightarrow$ (2).
\end{proof}

Lemma \ref{Lem1} and Theorem \ref{Th1} show that if $T$ and $\nabla$ are a
section of $\tau_{0}^{\dag}$ and a flat connection in $\pi_{l}$, respectively,
satisfying condition (\ref{HJCond}), then $\nabla$ is completely determined by
$T$ via $\nabla=\nabla^{T}$. This remark motivates the following.

\begin{definition}
A section $T$ of $\tau_{0}^{\dag}$ is a solution of Problem \ref{HJProb} iff
$\nabla^{T}$ is flat and $(T,\nabla^{T})$ satisfies condition (\ref{HJCond}).
\end{definition}

If a solution $T$ of Problem \ref{HJProb} exists, then $J^{l}$ is (locally)
foliated by $l$th jet prolongations of solutions of the Euler-Lagrange
equations. Moreover, such prolongations are projections of solutions of the
HDWE whose image lies in $\operatorname{im}T$.

Notice that assertions (1) and (2) in Theorem \ref{Th1} are still equivalent if
$Z\subset J^{1}\tau^{\dag}$ is the geometric portrait of a generic system
$\mathscr{E}$ of first order PDEs, i.e., a subbundle of $(\tau^{\dag})_{1,0}$.
On the other hand, when $\mathscr{E}$ are HDWE, solutions of Problem \ref{HJProb} can be
explicitly characterized in terms of the Hamiltonian section $\mathscr{H}$ as
follows. First of all, notice that $T_{\mathscr{H}}:=\mathscr{H}\circ
T:J^{l}\longrightarrow\mathscr{M}$ is a differential form in $\Lambda
_{n-1}^{n}(J^{l},\pi_{l})$. It holds the following

\begin{proposition}
A section $T:J^{l}\longrightarrow J^{\dag}$ of $\tau_{0}^{\dag}$ is a solution
of Problem \ref{HJProb} iff
\begin{equation}
i^{1,n}(\nabla^{T})(dT_{\mathscr{H}})=0, \label{HamJac1}%
\end{equation}
and $\nabla^{T}$ is flat.
\end{proposition}

\begin{proof}
Let $T$ be a section of $\tau_{0}^{\dag}$. Then $\operatorname{im}%
(T_{[1]}\circ\nabla^{T})\subset Z$ iff $i^{1,n}(T_{[1]}\circ\nabla^{T}%
)(\omega|_{T})=0$. Indeed, for $y\in J^{l}$,
\[
i^{1,n}(T_{[1]}\circ\nabla^{T})(\omega|_{T})(y)=i^{1,n}((T_{[1]}\circ
\nabla^{T})(y))(\omega|_{T(y)})=0
\]
iff $(T_{[1]}\circ\nabla^{T})(y)\in\operatorname{Ker}\omega|_{T(y)}\subset Z$.
Thus, $T$ is a solution of Problem \ref{HJProb} iff
\begin{align*}
0  &  =i^{1,n}(T_{[1]}\circ\nabla^{T})(\omega|_{T})\\
&  =i^{1,n}(\nabla^{T})(T^{\ast}(\omega))\\
&  =i^{1,n}(\nabla^{T})(T^{\ast}(d\theta))\\
&  =i^{1,n}(\nabla^{T})(d(T^{\ast}\circ\mathscr{H}^{\ast})(\Theta))\\
&  =i^{1,n}(\nabla^{T})(d(\mathscr{H}\circ T)^{\ast}(\Theta))\\
&  =i^{1,n}(\nabla^{T})(dT_{\mathscr{H}}).
\end{align*}
\end{proof}

Notice that, if $T$ is a section of $\tau_{0}^{\dag}$ but $\nabla=\nabla^{T}$
is not flat, then (\ref{HamJac1}) is still a sufficient condition for
$(T,\nabla)$ to satisfy condition (\ref{HJCond}).

\section{Hamilton-Jacobi Formalism\label{SecHJForm}}

In practice, it is quite hard to find solutions of the generalized HJ problem. That's why
it is convenient to formulate a simpler problem.

\begin{problem}
[Hamilton-Jacobi]\label{HJProb2}The HJ problem consists in finding solutions
$T$ of the generalized HJ problem such that $d^{V}T=0$.
\end{problem}

Problem \ref{HJProb2} extends the standard HJ problem in
Lagrangian-Hamiltonian mechanics to higher order Lagrangian field theory.
Clearly, in view of the universal property of the tautological element
$\underline{\Theta}\in V\!\Lambda^{1}(J^{\dag},\tau^{\dag})\otimes
\Lambda_{n-1}^{n-1}(J^{l},\pi_{l})$, condition $d^{V}T=0$ is equivalent to
$T^{\ast}(d^{V}\underline{\Theta})=0$, which, in its turn, is equivalent to
$T^{\ast}(\omega)\in\Lambda_{n}^{n}(J^{l},\pi_{l})$ (see, for instance,
\cite{v09b}). Thus, if $d^{V}T=0$, we have
\[
i^{1,n}(\nabla^{T})(dT_{\mathscr{H}})=i^{1,n}(\nabla^{T})(T^{\ast}%
(\omega))=T^{\ast}(\omega).
\]
Collecting results in the previous section with the above remarks we get

\begin{theorem}
[Hamilton-Jacobi]\label{HJTheor}A section $T:J^{l}\longrightarrow J^{\dag}$ of
$\tau_{0}^{\dag}$ is a solution of Problem \ref{HJProb2} iff
\begin{equation}
T^{\ast}(\omega)=0 \tag{$\star\star$}\label{HJEq}%
\end{equation}
and $\nabla^{T}$ is flat.
\end{theorem}

Theorem \ref{HJTheor} generalizes the standard HJ theorem in
Lagrangian-Hamiltonian mechanics to higher order Lagrangian field theory and
we will refer to (\ref{HJEq}) as the \emph{HJ equations }(see below). If $T$
is a section of $\tau_{0}^{\dag}$ but $\nabla=\nabla^{T}$ is not flat, then
(\ref{HJEq}) is still a sufficient condition for i) $(T,\nabla)$ to satisfy
condition (\ref{HJCond}) of Problem \ref{HJProb} and ii) $T$ to be $d^{V}\!$-closed.

Notice that the HJ formalism presented in this and the previous sections is
actually independent of the fact that the Hamiltonian section is determined by
a Lagrangian field theory and, therefore, remains valid when $\mathscr{H}$ is
any section of the bundle $\mathscr{M}\longrightarrow J^{\dag}$. However,
Hamiltonian sections coming from Lagrangian theories (see the appendix) are of
a very special kind and the corresponding HJ formalism is much simpler.

In the remaining part of this section we find coordinate formulations of both
Problems \ref{HJProb} and \ref{HJProb2}. Let $\ldots,x^{i},\ldots
,u_{I}^{\alpha},\ldots$, $|I|{}\leq l$ be jet coordinates on $J^{l}$ and
$\ldots,p_{\alpha}^{I.i},\ldots$ associated fiber coordinates on $J^{\dag}$.
Let $T$ be a section of $\tau_{0}^{\dag}$. Locally,
\[
T=%
%TCIMACRO{\tsum \nolimits_{|I|{}\leq l}}%
%BeginExpansion
{\textstyle\sum\nolimits_{|I|{}\leq l}}
%EndExpansion
T_{\alpha}^{I.i}d^{V}\!u_{I}^{\alpha}\otimes d^{n-1}x_{i}.
\]
The symbols of $\nabla^{T}$ are
\begin{align*}
\nabla_{I.i}^{\alpha}  &  :=(\nabla^{T})^{\ast}(u_{I}^{\alpha}{}_{|i})\\
&  =(\mathbb{F\!}\mathscr{H}\circ T)^{\ast}(u_{I}^{\alpha}{}_{|i})\\
&  =(T^{\ast}\circ\mathbb{F\!}\mathscr{H}^{\ast})(u_{I}^{\alpha}{}_{|i})\\
&  =\tfrac{\partial H}{\partial p_{\alpha}^{I.i}}\circ T\\
&  =\left\{
\begin{array}
[c]{ll}%
u_{Ii}^{\alpha} & \text{if }|I|{}<l\\
\mathfrak{s}_{Ii}^{\alpha}\circ T & \text{if }|I|{}=l
\end{array}
\right.  .
\end{align*}
For $|I|{}=l$, put $\ldots,\overline{\mathfrak{s}}_{Ii}^{\alpha}%
:=\mathfrak{s}_{Ii}^{\alpha}\circ T,\ldots$. The curvature of $\nabla^{T}$ is
then
\[
R^{T}:=%
%TCIMACRO{\tsum \nolimits_{|I|{}\leq l}}%
%BeginExpansion
{\textstyle\sum\nolimits_{|I|{}\leq l}}
%EndExpansion
R_{I.ij}^{\alpha}dx^{i}dx^{j}\otimes\tfrac{\partial}{\partial u_{I}^{\alpha}%
},
\]
with
\[
R_{I.ij}^{\alpha}=\tfrac{1}{2}(D_{i}\nabla_{I.j}^{\alpha}-D_{j}\nabla
_{I.i}^{\alpha})\circ\nabla^{T}=\left\{
\begin{array}
[c]{ll}%
0 & \text{if }|I|{}<l\\
\tfrac{1}{2}(D_{i}\overline{\mathfrak{s}}_{Ij}^{\alpha}-D_{j}\overline
{\mathfrak{s}}_{Ii}^{\alpha})\circ\nabla^{T} & \text{if }|I|{}=l
\end{array}
\right.  .
\]
Recall that the $\mathfrak{s}_{Ij}^{\alpha}$'s, $|I|{}=l$, are implicitly
defined by Eq.{} (\ref{ImpEq}) and therefore
\[%
%TCIMACRO{\tsum \nolimits_{|I|{}=l}}%
%BeginExpansion
{\textstyle\sum\nolimits_{|I|{}=l}}
%EndExpansion
(\partial_{\alpha}^{Kk}{}\partial_{\beta}^{Ii}L\,D_{j}\overline{\mathfrak{s}%
}_{Ii}^{\alpha})\circ\nabla^{T}=\delta_{Jl}^{Kk}D_{j}T_{\beta}^{J.l}%
\circ\nabla^{T}.
\]
Let $\mathbf{M}=\left\Vert M_{I}^{\alpha}{}_{J}^{\beta}\right\Vert^{(\alpha,I)}_{(\beta,J)}$ be the inverse matrix of
$\mathbf{H}$, $|I|,|J|{}=l+1$. Then
\[
(D_{j}\overline{\mathfrak{s}}_{Ii}^{\alpha})\circ\nabla^{T}=%
%TCIMACRO{\tsum \nolimits_{|J|{}=l}}%
%BeginExpansion
{\textstyle\sum\nolimits_{|J|{}=l}}
%EndExpansion
(M_{Jl}^{\beta}{}{}_{Ii}^{\alpha}D_{j}T_{\beta}^{J.l})\circ\nabla^{T}, \quad |I|{}=l.
\]
Therefore
\[
R_{I.ij}^{\alpha}=\left\{
\begin{array}
[c]{ll}%
0 & \text{if }|I|{}<l\\%
%TCIMACRO{\tsum \nolimits_{|J|{}=l}}%
%BeginExpansion
{\textstyle\sum\nolimits_{|J|{}=l}}
%EndExpansion
(M_{Jl}^{\beta}{}{}_{I[i}^{\alpha}D_{j]}T_{\beta}^{J.l})\circ\nabla^{T} &
\text{if }|I|{}=l
\end{array}
\right.  .
\]
In particular, $\nabla^{T}$ is flat iff
\[%
%TCIMACRO{\tsum \nolimits_{|J|{}=l}}%
%BeginExpansion
{\textstyle\sum\nolimits_{|J|{}=l}}
%EndExpansion
(M_{Jl}^{\beta}{}{}_{I[i}^{\alpha}D_{j]}T_{\beta}^{J.l})\circ\nabla
^{T}=0,\quad|I|{}=l,
\]
where square brackets denote skew-symmetrization.

Now
\[
T_{\mathscr{H}}=%
%TCIMACRO{\tsum \nolimits_{|I|{}\leq l}}%
%BeginExpansion
{\textstyle\sum\nolimits_{|I|{}\leq l}}
%EndExpansion
T_{\alpha}^{I.i}du_{I}^{\alpha}d^{n-1}x_{i}-(H\circ T)d^{n}x
\]
and
\begin{align*}
&  dT_{\mathscr{H}}\\
&  =%
%TCIMACRO{\tsum \nolimits_{|I|{}\leq l}}%
%BeginExpansion
{\textstyle\sum\nolimits_{|I|{}\leq l}}
%EndExpansion
dT_{\alpha}^{I.i}du_{I}^{\alpha}d^{n-1}x_{i}-d(H\circ T)d^{n}x\\
&  =%
%TCIMACRO{\tsum \nolimits_{|I|,|J|{}\leq l}}%
%BeginExpansion
{\textstyle\sum\nolimits_{|I|,|J|{}\leq l}}
%EndExpansion
\partial_{\beta}^{J}T_{\alpha}^{I.i}du_{J}^{\beta}du_{I}^{\alpha}d^{n-1}%
x_{i}\\
&  -%
%TCIMACRO{\tsum \nolimits_{|J|{}\leq l}}%
%BeginExpansion
{\textstyle\sum\nolimits_{|J|{}\leq l}}
%EndExpansion
[\partial_{i}T_{\beta}^{J.i}+(\partial_{\beta}^{J}H)\circ T+%
%TCIMACRO{\tsum \nolimits_{|I|{}\leq l}}%
%BeginExpansion
{\textstyle\sum\nolimits_{|I|{}\leq l}}
%EndExpansion
(\tfrac{\partial H}{\partial p_{\alpha}^{I.i}}\circ T)\partial_{\beta}%
^{J}T_{\alpha}^{I.i}]du_{J}^{\beta}d^{n}x,
\end{align*}
so that
\begin{align*}
&  i^{1,n}(\nabla^{T})(dT_{\mathscr{H}})\\
&  =%
%TCIMACRO{\tsum \nolimits_{|J|{}\leq l}}%
%BeginExpansion
{\textstyle\sum\nolimits_{|J|{}\leq l}}
%EndExpansion
[%
%TCIMACRO{\tsum \nolimits_{|I|{}\leq l}}%
%BeginExpansion
{\textstyle\sum\nolimits_{|I|{}\leq l}}
%EndExpansion
\overline{\mathfrak{s}}_{Ii}^{\alpha}(\partial_{\beta}^{J}T_{\alpha}%
^{I.i}-\partial_{\alpha}^{I}T_{\beta}^{J.i})-\partial_{i}T_{\beta}%
^{J.i}-(\partial_{\beta}^{J}H)\circ T-%
%TCIMACRO{\tsum \nolimits_{|I|{}\leq l}}%
%BeginExpansion
{\textstyle\sum\nolimits_{|I|{}\leq l}}
%EndExpansion
(\tfrac{\partial H}{\partial p_{\alpha}^{I.i}}\circ T)\partial_{\beta}%
^{J}T_{\alpha}^{I.i}]du_{J}^{\beta}d^{n}x\\
&  =-%
%TCIMACRO{\tsum \nolimits_{|J|{}\leq l}}%
%BeginExpansion
{\textstyle\sum\nolimits_{|J|{}\leq l}}
%EndExpansion
[D_{i}T_{\beta}^{J.i}\circ\nabla^{T}+(\partial_{\beta}^{J}H)\circ
T]du_{J}^{\beta}d^{n}x.
\end{align*}
Thus $T$ is a solution of Problem \ref{HJProb} iff
\begin{align*}%
%TCIMACRO{\tsum \nolimits_{|J|{}=l}}%
%BeginExpansion
{\textstyle\sum\nolimits_{|J|{}=l}}
%EndExpansion
(M_{Jl}^{\beta}{}{}_{I[i}^{\alpha}D_{j]}T_{\beta}^{J.l})\circ\nabla^{T}  &
=0,\quad|I|{}=l,\\
D_{i}T_{\beta}^{J.i}\circ\nabla^{T}+(\partial_{\beta}^{J}H)\circ T  &
=0,\quad|J|{}\leq l.
\end{align*}
Now, suppose that $d^{V}\!T=0$. Then, locally, $T=d^{V}\!S$ for some $S\in
\Lambda_{n-1}^{n-1}(J^{l},\pi_{l})$. If $S$ is locally given by
\[
S=S^{i}d^{n-1}x_{i},
\]
then, locally,
\[
T=\partial_{\beta}^{J}S^{i}d^{V}\!u_{J}^{\beta}\otimes d^{n-1}x_{i},
\]
and $T^{\ast}(\omega)$ is locally given by
\[
T^{\ast}(\omega)=\partial_{\beta}^{J}(\partial_{i}S^{i}+H\circ d^{V}\!
S)du_{J}^{\beta}d^{n}x.
\]
Therefore, the HJ equations read
\[
\partial_{\beta}^{J}(\partial_{i}S^{i}+H\circ d^{V}\!S)=0,\quad|J|{}\leq l
\]
and $T$ is a solution of Problem \ref{HJProb2} iff
\begin{align*}%
%TCIMACRO{\tsum \nolimits_{|J|{}=l}}%
%BeginExpansion
{\textstyle\sum\nolimits_{|J|{}=l}}
%EndExpansion
(M_{Jl}^{\beta}{}{}_{I[i}^{\alpha}D_{j]}T_{\beta}^{J.l})\circ\nabla^{T}  &
=0,\quad|I|{}=l,\\
\partial_{\beta}^{J}(\partial_{i}S^{i}+H\circ d^{V}\!S)  &  =0,\quad|J|{}\leq l.
\end{align*}

\section{The HJ Formalisms of Equivalent Lagrangian
Theories\label{SecEquivLag}}

We recall that two Lagrangian densities differing by a total divergence
determine the same Euler-Lagrange equations. We may then wonder how do Problems
\ref{HJProb} and \ref{HJProb2}, and their solutions, change when adding a
total divergence to the Lagrangian density. In order to answer this question
let us first discuss how does the Lagrangian-Hamiltonian picture depicted in
Sections \ref{SecLHForm} and \ref{SecInvLeg} changes when adding a total
divergence to the Lagrangian density (see also \cite{v09}).

Let $(\pi,\mathscr{L})$ be a Lagrangian field theory. Recall the geometric
definition of total divergence. Let $\varrho\in\Lambda_{n-1}^{n-1}(J^{l}%
,\pi_{l})$, and consider $d\varrho\in\Lambda^{n}(J^{l})$. We put
\[
\overline{d}\varrho:=i^{0,n}(\mathscr{C})(d\varrho)\in\Lambda_{n}^{n}%
(J^{l+1},\pi_{l}).
\]
$\overline{d}\varrho$ is then the \emph{total divergence }of $\varrho$. If
$\varrho$ is locally given by $\varrho=\varrho^{i}d^{n-1}x_{i}$, $\overline{d}\varrho$ is
locally given by $\overline{d}\varrho=D_{i}\varrho^{i}d^{n}x$. Then the Lagrangian field
theory $(\pi,\tilde{\mathscr{L}})$, where $\tilde{\mathscr{L}}%
:=\mathscr{L}+\overline{d}\varrho$, determines the same Euler-Lagrange equations as
$(\pi,\mathscr{L})$. For other geometric objects determined by $\tilde
{\mathscr{L}}$ (Hamiltonian section, inverse Legendre transform, etc.) we use
the same notation as for those determined by $\mathscr{L}$ simply adding a
tilde. For instance, we denote by $\tilde{\mathscr{H}}$ the Hamiltonian section
determined by $\tilde{\mathscr{L}}$ (see below).

Now, $d^{V}\!\varrho$ is a section of $\tau_{0}^{\dag}$. It determines an
automorphism $\Psi:J^{\dag}\longrightarrow J^{\dag}$ of $\tau_{0}^{\dag}$ via
\[
\Psi(P):=P-d^{V}\!\varrho|_{y},\quad P\in J^{\dag},\;y:=\tau_{0}^{\dag}(P)\in
J^{l}.
\]
Similarly, the section $\pi_{l+1,l}^{\circ}(d^{V}\!\varrho)$ determines an
automorphism $\Psi^{\prime}:\pi_{l+1,l}^{\circ}(J^{\dag})\longrightarrow
\pi_{l+1,l}^{\circ}(J^{\dag})$ of $\mathfrak{p}$. Notice that diagram
\[%
%TCIMACRO{\TeXButton{TeX field}{\xymatrix{\pi_{l+1,l}^\circ( J^\dag) \ar
%[d] \ar[r]^-{\Psi^\prime} & \pi_{l+1,l}^\circ( J^\dag) \ar[d] \\
%J^\dag\ar[r]_-{\Psi} & J^\dag}}}%
%BeginExpansion
\xymatrix{\pi_{l+1,l}^\circ( J^\dag) \ar[d] \ar[r]^-{\Psi^\prime}
& \pi_{l+1,l}^\circ( J^\dag) \ar[d] \\
J^\dag\ar[r]_-{\Psi} & J^\dag}%
%EndExpansion
\]
commutes, i.e., the pair $(\Psi,\Psi^{\prime})$ is also a morphism of bundles
over $J^{\dag}$.

\begin{theorem}
\quad i) $\tilde{\theta}^{\prime}=\Psi^{\prime}{}^{\ast}(\theta^{\prime}%
)+(\pi_{l+1,l}\circ\mathfrak{p})^{\ast}(d\varrho)$, ii) $\tilde{\omega
}^{\prime}=\Psi^{\prime}{}^{\ast}(\omega^{\prime})$ and iii) $\tilde
{\mathscr{P}}=\Psi^{\prime}{}^{-1}(\mathscr{P})$. Moreover, $(\pi
,\mathscr{L})$ is hyperregular iff $(\pi,\tilde{\mathscr{L}})$ is
hyperregular, and, in this case, iv) $\tilde{\mathfrak{s}}=\Psi^{\prime}%
{}^{-1}\circ\mathfrak{s}\circ\Psi$, v) $\mathsf{\tilde{\mathsf{leg}}}{}%
^{-1}=\mathsf{leg}^{-1}\circ\Psi$, vi) $\tilde{\theta}=\Psi(\theta)+\tau
_{0}^{\dag}{}^{\ast}(d\varrho)$, vii) $\tilde{\omega}=\Psi(\omega)$.
\end{theorem}

\begin{proof}
First of all compute
\begin{align*}
\Psi^{\prime}{}^{\ast}(\theta^{\prime})  &  =(\Psi^{\prime}{}^{\ast}\circ
h^{\ast})(\Theta)+(\Psi^{\prime}{}^{\ast}\circ\mathfrak{p}{}^{\ast
})(\mathscr{L})\\
&  =(h\circ\Psi^{\prime})^{\ast}(\Theta)+(\mathfrak{p}\circ\Psi^{\prime}%
{})^{\ast}(\mathscr{L})\\
&  =(h\circ\Psi^{\prime})^{\ast}(\Theta)+\mathfrak{p}{}^{\ast}(\mathscr{L}).
\end{align*}
Now, locally
\begin{align*}
(h\circ\Psi^{\prime})^{\ast}(p_{\alpha}^{I.i})  &  =p_{\alpha}^{I.i}%
-\partial_{\alpha}^{I}\varrho^{i},\quad|I|{}\leq l\\
(h\circ\Psi^{\prime})^{\ast}(p)  &  =%
%TCIMACRO{\tsum \nolimits_{|I|{}\leq l}}%
%BeginExpansion
{\textstyle\sum\nolimits_{|I|{}\leq l}}
%EndExpansion
(p_{\alpha}^{I.i}-\partial_{\alpha}^{I}\varrho^{i})u_{Ii}^{\alpha}.
\end{align*}
Thus,
\begin{align*}
(h\circ\Psi^{\prime})^{\ast}(\Theta)  &  =%
%TCIMACRO{\tsum \nolimits_{|I|{}\leq l}}%
%BeginExpansion
{\textstyle\sum\nolimits_{|I|{}\leq l}}
%EndExpansion
(p_{\alpha}^{I.i}-\partial_{\alpha}^{I}\varrho^{i})du_{I}^{\alpha}d^{n-1}%
x_{i}-%
%TCIMACRO{\tsum \nolimits_{|I|{}\leq l}}%
%BeginExpansion
{\textstyle\sum\nolimits_{|I|{}\leq l}}
%EndExpansion
(p_{\alpha}^{I.i}-\partial_{\alpha}^{I}\varrho^{i})u_{Ii}^{\alpha}d^{n}x\\
&  =h^{\ast}(\Theta)-\mathfrak{p}{}^{\ast}(\pi_{l+1,l}^{\ast}(d\varrho
)-\overline{d}\varrho)\\
&  =h^{\ast}(\Theta)-(\pi_{l+1,l}\circ\mathfrak{p})^{\ast}(d\varrho
)+\mathfrak{p}{}^{\ast}(\overline{d}\varrho)
\end{align*}
and, therefore,
\begin{align*}
\Psi^{\prime}{}^{\ast}(\theta^{\prime})  &  =h^{\ast}(\Theta)-(\pi
_{l+1,l}\circ\mathfrak{p})^{\ast}(d\varrho)+\mathfrak{p}{}^{\ast}(\overline
{d}\varrho)+\mathfrak{p}{}^{\ast}(\mathscr{L})\\
&  =h^{\ast}(\Theta)+\mathfrak{p}{}^{\ast}(\tilde{\mathscr{L}})-(\pi
_{l+1,l}\circ\mathfrak{p})^{\ast}(d\varrho)\\
&  =\tilde{\theta}^{\prime}-(\pi_{l+1,l}\circ\mathfrak{p})^{\ast}(d\varrho).
\end{align*}
We then have $\tilde{\omega}^{\prime}=d\tilde{\theta}^{\prime}=\Psi^{\prime}%
{}^{\ast}(\omega^{\prime})$. It follows that $\tilde{\mathscr{P}}=\Psi
^{\prime}{}^{-1}(\mathscr{P})$. Indeed, let $z\in\pi_{l+1,l}^{\circ}(J^{\dag
})$ and $\Pi$ be an element over $z$ of the first jet bundle of $\pi
_{l+1,l}^{\circ}(J^{\dag})\longrightarrow M$. Then
\[
i^{1,n}(\Pi)(\tilde{\omega}^{\prime}|_{z})=i^{1,n}(\Pi)(\Psi^{\prime}{}^{\ast
}(\omega^{\prime})|_{z})=i^{1,n}(\Psi_{\lbrack1]}^{\prime}(\Pi))(\omega
^{\prime}|_{\Psi^{\prime}(z)}),
\]
so that $z\in\tilde{\mathscr{P}}$ iff $\Psi^{\prime}(z)\in\mathscr{P}$. Since
$(\Psi,\Psi^{\prime})$ is an isomorphism of bundles over $J^{\dag}$,
$\tilde{\mathscr{P}}$ projects diffeomorphically to $J^{\dag}$ iff
$\mathscr{P}$ does, i.e., $(\pi,\mathscr{L})$ is hyperregular iff $(\pi
,\tilde{\mathscr{L}})$ is. In this case, it immediately follows from iii) that
$\tilde{\mathfrak{s}}=\Psi^{\prime}{}^{-1}\circ\mathfrak{s}\circ\Psi$. Then
\begin{align*}
\mathsf{\tilde{\mathsf{leg}}}{}^{-1}  &  =\mathfrak{p}\circ\tilde
{\mathfrak{s}}\\
&  =\mathfrak{p}\circ\Psi^{\prime}{}^{-1}\circ\mathfrak{s}\circ\Psi\\
&  =\mathfrak{p}\circ\mathfrak{s}\circ\Psi\\
&  =\mathsf{leg}^{-1}\circ\Psi,
\end{align*}
and
\begin{align*}
\tilde{\theta}  &  =\tilde{\mathfrak{s}}^{\ast}(\tilde{\theta}^{\prime})\\
&  =(\Psi^{\prime}{}^{-1}\circ\mathfrak{s}\circ\Psi)^{\ast}(\Psi^{\prime}%
{}^{\ast}(\theta^{\prime})+(\pi_{l+1,l}\circ\mathfrak{p})^{\ast}(d\varrho))\\
&  =(\Psi^{\prime}{}\circ\Psi^{\prime}{}^{-1}\circ\mathfrak{s}\circ\Psi
)^{\ast}(\theta^{\prime})+(\pi_{l+1,l}\circ\mathfrak{p}\circ\Psi^{\prime}%
{}^{-1}\circ\mathfrak{s}\circ\Psi)^{\ast}(d\varrho)\\
&  =(\Psi^{\ast}\circ\mathfrak{s}^{\ast})(\theta^{\prime})+(\pi_{l+1,l}%
\circ\mathfrak{p}\circ\mathfrak{s}\circ\Psi)^{\ast}(d\varrho)\\
&  =\Psi^{\ast}(\theta)+(\pi_{l+1,l}\circ\mathsf{leg}^{-1}\circ\Psi)^{\ast
}(d\varrho)\\
&  =\Psi^{\ast}(\theta)+(\tau_{0}^{\dag}\circ\Psi)^{\ast}(d\varrho)\\
&  =\Psi^{\ast}(\theta)+\tau_{0}^{\dag}{}^{\ast}(d\varrho).
\end{align*}
We then have $\tilde{\omega}=d\tilde{\theta}=\Psi^{\ast}(\omega)$.
\end{proof}

Now, let $\tilde{T}$ be a section of $\tau_{0}^{\dag}$. We put $\tilde{\nabla
}^{\tilde{T}}:=$ $\mathsf{\tilde{\mathsf{leg}}}{}^{-1}\circ\tilde{T}$. We then
have the following

\begin{corollary}
Let $T$ be a section of $\tau_{0}^{\dag}$ and $\tilde{T}:=T+d^{V}\!\varrho$.
$T$ is a solution of the (generalized) HJ Problem determined by $\mathscr{L}$
iff $\tilde{T}$ is a solution of the (generalized) HJ Problem determined by
$\tilde{\mathscr{L}}$.
\end{corollary}

\begin{proof}
Clearly, $\tilde{T}=\Psi^{-1}\circ T$. Therefore,
\begin{align*}
\tilde{\nabla}^{\tilde{T}}  &  =\mathsf{\tilde{\mathsf{leg}}}{}^{-1}%
\circ\tilde{T}\\
&  =\mathsf{\tilde{\mathsf{leg}}}{}^{-1}\circ\Psi^{-1}\circ T\\
&  =\mathsf{\mathsf{leg}}^{-1}\circ T\\
&  =\nabla^{T},
\end{align*}
and
\begin{align*}
d\tilde{T}_{\tilde{\mathscr{H}}}  &  =\tilde{T}^{\ast}(\tilde{\omega})\\
&  =((\Psi^{-1}\circ T)^{\ast}\circ\Psi^{\ast})(\omega)\\
&  =T^{\ast}(\omega)\\
&  =dT_{\mathscr{H}},
\end{align*}
so that
\[
i^{1,n}(\tilde{\nabla}^{\tilde{T}})(d\tilde{T}_{\tilde{\mathscr{H}}}%
)=i^{1,n}(\nabla^{T})(dT_{\mathscr{H}}).
\]
Finally, $d^V\! T=0$ iff $d^V \! \tilde{T}=0$.
\end{proof}

We conclude that equivalent Lagrangians (of the same order) determine
essentially equivalent (generalized) HJ Problems.

\section{An Example: the Biharmonic Equation\label{SecExamp}}

Consider the second order Lagrangian field theory $(\pi,\mathscr{L})$, where
$\pi$ is the bundle $\mathbb{R}^{n}\times\mathbb{R}\ni(\ldots,x^{i}%
,\ldots;u)\longmapsto(\ldots,x^{i},\ldots)\in\mathbb{R}^{n}$ and
$\mathscr{L}=Ld^{n}x$, with
\[
L=\tfrac{1}{2}u_{ij}u^{ij}.
\]
Throughout this section indexes are lowered and raised using the standard
Euclidean metric on $\mathbb{R}^{n}$. The Euler-Lagrange equation reads
\begin{equation}
\nabla^{4}u=0, \label{BHE}%
\end{equation}
where $\nabla^{4}$ is the bilaplacian:
\[
\nabla^{4}u:=u,_{ij}{}^{ij}.
\]
Equation (\ref{BHE}) is known as the \emph{biharmonic equation} and its
solutions as \emph{biharmonic functions}. For $n=2$ the biharmonic equation is
often used to model the bending of a thin, elastic, unclamped plate subjected
to given boundary conditions. In this case $u$ represents the (small)
displacement from the plane (whose points are labeled by Cartesian coordinates
$x^{1},x^{2}$) where the plate is and $L$ is the \emph{bending energy} of the plate.

Denote by $\ldots,p^{.i},\ldots,p^{j.i},\ldots$ natural coordinates in
$J^{\dag}\pi_{1}$, and by $p$ the remaining coordinate in $\mathscr{M}\pi_{1}%
$. $\mathscr{P}$ is described by equations
\begin{equation}
u_{ij}=p^{(i.j)}, \label{P}%
\end{equation}
where round brackets denotes symmetrization, namely
$p^{(i.j)}:=\tfrac{1}{2}(p^{i.j}+p^{j.i})$. Equations (\ref{P}) show that the
theory is hyperregular. The Hamiltonian section is given by $\mathscr{H}^{\ast
}(p)=H$ with
\[
H=p^{.i}u_{i}+\tfrac{1}{2}p_{(i.j)}p^{(i.j)},
\]
and the HDWE read
\begin{equation}
\left\{
\begin{array}
[c]{l}%
p^{.i},_{i}=0\\
p^{j.i},_{i}=p^{.j}\\
u,_{i}=u_{i}\\
u_{i},_{j}=p_{(i.j)}%
\end{array}
\right.  ,\quad i,j=1,\ldots,n. \label{HDWExamp}%
\end{equation}
which clearly cover (\ref{BHE}).

The HJ equation reads
\begin{equation}
\partial_{i}S^{i}+u_{i}\tfrac{\partial}{\partial u}S^{i}+\tfrac{1}{4}%
(\tfrac{\partial}{\partial u^{i}}S_{j}+\tfrac{\partial}{\partial u^{j}}%
S_{i})\tfrac{\partial}{\partial u_{i}}S^{j}=f(x). \label{HJExamp}%
\end{equation}
Now, let $\phi=\phi(x)$ be a biharmonic function, i.e., $\nabla^{4}\phi=0$.
Then a direct computation shows that
\[
S^{i}:=u^{j}\phi,_{j}{}^{i}-u\phi,_{j}{}^{ji}+G^{i}(x),\quad i=1,\ldots,n
\]
is a solution of (\ref{HJExamp}), for every $n$tuple of functions
$\ldots,G^{i}=G^{i}(x),\ldots$ of the $x$'s only (and we can even achieve
$f=0$ by putting $\ldots,G^{i}:=\tfrac{1}{2}(\phi\phi,_{j}{}^{ji}-\phi^{j}%
\phi,_{j}{}^{i}),\ldots$). Now
\[
T:=d^{V}\!(S^{i}d^{n-1}x_{i})=(-\phi,_{j}{}^{ji}d^{V}\!u+\phi,{}^{ij}%
d^{V}\!u_{j})\otimes d^{n-1}x_{i},
\]
and the symbols of $\nabla^{T}$ are
\begin{align*}
(\nabla^{T})^{\ast}(u_{i})  &  =u_{i},\\
(\nabla^{T})^{\ast}(u_{i|j})  &  =(\nabla^{T})^{\ast}(u_{j|i})=\phi,_{ij}.
\end{align*}
In particular, $\nabla^{T}$ is flat and $J^{1}\pi$ is foliated by $\nabla^{T}%
$-flat sections:
\begin{align}
u  &  =\phi+A_{i}x^{i}+B,\nonumber\\
u_{i}  &  =\phi,_{i}+A_{i}, \label{FS}%
\end{align}
where $\ldots,A_{i},\ldots,B$ are constants. $\nabla^{T}$-constant sections
(\ref{FS}) are first jet prolongations of the biharmonic functions
\[
\phi_{A,B}:=\phi+A_{i}x^{i}+B.
\]
Finally for all $\ldots,A_{i},\ldots,B$, $T\circ\phi_{A,B}$ is given by
\begin{align*}
p^{.i}  &  =-\phi,_{j}{}^{ji}\\
p^{j.i}  &  =\phi,{}^{ij}\\
u  &  =\phi+A_{i}x^{i}+B\\
u_{i}  &  =\phi,_{i}+A_{i}%
\end{align*}
and is a solution of (\ref{HDWExamp}).

\appendix{}

\section{Hamiltonian Theories with Lagrangian Counterparts}

A section $\mathscr{H}$ of the bundle $\mathscr{M}\longrightarrow J^{\dag}$
determines a PD-Hamiltonian system $d\theta$, $\theta:=\mathscr{H}^{\ast}(\Theta)$. It is
natural to understand $\theta$ as an (higher order) Hamiltonian field theory.
However, as noticed in Section \ref{SecHJForm}, not all $\mathscr{H}$'s come
from a Lagrangian field theory, i.e., not all Hamiltonian field theories are
determined by an underlying Lagrangian field theory. In this appendix we
characterize Hamiltonian field theories having a hyperregular Lagrangian
counterpart. Notice, preliminarily, that any section $\mathscr{H}$ of
$\mathscr{M}\longrightarrow J^{\dag}$ determines a morphism of bundles over
$J^{l}$ (fiber derivative \cite{fr05}): $\mathbb{F}\!\mathscr{H}:J^{\dag
}\longrightarrow J^{1}\pi_{l}$. Recall that, if $\mathscr{M}$ is locally given
by $\mathscr{H}^{\ast}(p)=H$, then $\mathbb{F}\!\mathscr{H}$ is given by
$\ldots,\mathbb{F}\!\mathscr{H}(u_{I}^{\alpha}{}_{|i})=\tfrac{\partial
H}{\partial p_{\alpha}^{I.i}},\ldots$.

\begin{theorem}
Let $\mathscr{H}$ be a section of $\mathscr{M}\longrightarrow J^{\dag}$.
$\mathscr{H}$ is the Hamiltonian section of an hyperregular Lagrangian field
theory of the order $l+1$ iff $\mathbb{F}\!\mathscr{H}$ takes its values in
$J^{l+1}$ and $\mathbb{F}\!\mathscr{H}:J^{\dag}\longrightarrow J^{l+1}$ is a
surjective submersion with connected fibers.
\end{theorem}

\begin{proof}
We need to prove just the \textquotedblleft only if\textquotedblright%
\ implication. Thus, let $\mathbb{F}\!\mathscr{H}$ be a surjective submersion
onto $J^{l+1}$ with connected fibers. In particular, $\mathbb{F}\!\mathscr{H}$ uniquely
determines a section $\mathfrak{s}$ of $\pi_{l+1,l}^{\circ}(J^{\dag
})\longrightarrow J^{\dag}$ such that $\mathfrak{p}\circ\mathfrak{s}%
=\mathbb{F}\!\mathscr{H}$. Moreover,
\[
\tfrac{\partial H}{\partial p_{\alpha}^{I.i}}=u_{Ii}^{\alpha},\quad|I|{}<l
\]
and
\[
\tfrac{\partial H}{\partial p_{\alpha}^{I.i}}=\tfrac{\partial H}{\partial
p_{\alpha}^{J.j}},\quad\text{whenever }|I|{}=|J|{}=l\text{ and }Ii=Jj.
\]
Finally,
\begin{equation}
\operatorname{rank}\left\Vert
%TCIMACRO{\tsum \nolimits_{|J_{1}|{}=l}}%
%BeginExpansion
{\textstyle\sum\nolimits_{|J_{1}|{}=l}}
%EndExpansion
\tfrac{\partial^{2}H}{\partial p_{\alpha}^{I.i}\partial p_{\beta}^{J_{1}.j}%
}\delta_{J}^{J_{1}j}\right\Vert_{(J,\beta)}^{(I.i,\alpha)}=\dim\pi_{l+1,l}.
\label{HessH}%
\end{equation}
Now, put $\theta:=\mathscr{H}^{\ast}(\Theta)$, as above, and denote by
$\Delta\in V\!D(J^{\dag},\tau_{0}^{\dag})$ the Liouville vector field of
$J^{\dag}$. $\Delta$ is locally given by%
\[
\Delta=%
%TCIMACRO{\tsum \nolimits_{|I|{}\leq l}}%
%BeginExpansion
{\textstyle\sum\nolimits_{|I|{}\leq l}}
%EndExpansion
p_{\alpha}^{I.i}\tfrac{\partial}{\partial p_{\alpha}^{I.i}}.
\]
Notice that the differential form $\theta-L_{\Delta}\theta$ is actually in
$\Lambda_{n}^{n}(J^{\dag},\tau^{\dag})$. Indeed, locally
\begin{align*}
L_{\Delta}\theta &  =L_{\Delta}(%
%TCIMACRO{\tsum \nolimits_{|I|{}\leq l}}%
%BeginExpansion
{\textstyle\sum\nolimits_{|I|{}\leq l}}
%EndExpansion
p_{\alpha}^{I.i}du_{I}^{\alpha}d^{n-1}x_{i}-Hd^{n}x)\\
&  =%
%TCIMACRO{\tsum \nolimits_{|I|{}\leq l}}%
%BeginExpansion
{\textstyle\sum\nolimits_{|I|{}\leq l}}
%EndExpansion
\Delta(p_{\alpha}^{I.i})du_{I}^{\alpha}d^{n-1}x_{i}-\Delta(H)d^{n}x\\
&  =%
%TCIMACRO{\tsum \nolimits_{|I|{}\leq l}}%
%BeginExpansion
{\textstyle\sum\nolimits_{|I|{}\leq l}}
%EndExpansion
p_{\alpha}^{I.i}du_{I}^{\alpha}d^{n-1}x_{i}-\Delta(H)d^{n}x,
\end{align*}
and therefore, locally,
\[
\theta-L_{\Delta}\theta=(\Delta(H)-H)d^{n}x=(%
%TCIMACRO{\tsum \nolimits_{|I|{}\leq l}}%
%BeginExpansion
{\textstyle\sum\nolimits_{|I|{}\leq l}}
%EndExpansion
p_{\alpha}^{I.i}\tfrac{\partial H}{\partial p_{\alpha}^{I.i}}-H)d^{n}x.
\]
We claim that $\theta-L_{\Delta}\theta$ is the pull-back of a unique form
$\mathscr{L}\in\Lambda_{n}^{n}(J^{l+1},\pi_{l+1})$ via $\mathbb{F}\!\mathscr{H}$. To show
this, it is enough that $\theta-L_{\Delta}\theta$ is locally constant along
fibers of $\mathbb{F}\!\mathscr{H}$, i.e., $L_{Y}(\theta-L_{\Delta}\theta)=0$
for all $Y\in V\!D(J^{\dag},\mathbb{F}\!\mathscr{H})$. Let $Y$ be a vector field
over $J^{\dag}$ vertical with respect to $\mathbb{F}\!\mathscr{H}$. Then, in
particular, $Y$ is locally given by
\[
Y=%
%TCIMACRO{\tsum \nolimits_{|I|{}\leq l}}%
%BeginExpansion
{\textstyle\sum\nolimits_{|I|{}\leq l}}
%EndExpansion
Y_{\alpha}^{I.i}\tfrac{\partial}{\partial p_{\alpha}^{I.i}}%
\]
and, moreover,
\[
Y(\tfrac{\partial H}{\partial p_{\alpha}^{I.i}})=0.
\]
Let $Y$ be such a vector field. Locally,
\[
L_{Y}(\theta-L_{\Delta}\theta)=Y(%
%TCIMACRO{\tsum \nolimits_{|I|{}\leq l}}%
%BeginExpansion
{\textstyle\sum\nolimits_{|I|{}\leq l}}
%EndExpansion
p_{\alpha}^{I.i}\tfrac{\partial H}{\partial p_{\alpha}^{I.i}}-H)d^{n}x=(%
%TCIMACRO{\tsum \nolimits_{|I|{}\leq l}}%
%BeginExpansion
{\textstyle\sum\nolimits_{|I|{}\leq l}}
%EndExpansion
Y(p_{\alpha}^{I.i})\tfrac{\partial H}{\partial p_{\alpha}^{I.i}}%
-Y(H))d^{n}x=0.
\]

Now, let $\mathscr{L}\in\Lambda_{n}^{n}(J^{l+1},\pi_{l+1})$ be such that $\mathbb{F}\!%
\mathscr{H}^{\ast}(\mathscr{L})=\theta-L_{\Delta}\theta$. $\mathscr{L}$ is a
Lagrangian density of the order $l+1$. Prove that $(\pi,\mathscr{L})$ is a
hyperregular Lagrangian field theory whose Hamiltonian section is
$\mathscr{H}$. Let $\theta^{\prime}$ be the PD-Hamiltonian system on
$\pi_{l+1,l}^{\circ}(J^{\dag})\longrightarrow M$ determined by $\mathscr{L}$
and $\mathscr{P}$ its constraint submanifold. Then we have to prove that
$\mathscr{P}=\operatorname{im}\mathfrak{s}$. Since $\mathscr{P}$ and
$\operatorname{im}\mathfrak{s}$ are (closed) submanifolds (in $\pi
_{l+1,l}^{\circ}(J^{\dag})$) of the same dimension, it is enough to prove that
$\operatorname{im}\mathfrak{s}\subset\mathscr{P}$ in any local coordinate
neighborhood. If $\mathscr{L}$ is locally given by $\mathscr{L}=Ld^{n}x$,
then $\mathscr{P}$ is locally defined by (\ref{EqP}). Thus, $\operatorname{im}%
\mathfrak{s}\subset\mathscr{P}$ locally, iff
\[
\mathfrak{s}^{\ast}(\partial_{\alpha}^{I}L-%
%TCIMACRO{\tsum \nolimits_{|J|{}\leq l}}%
%BeginExpansion
{\textstyle\sum\nolimits_{|J|{}\leq l}}
%EndExpansion
\delta_{Ji}^{I}p_{\alpha}^{J.i})=0,\quad|I|{}=l+1.
\]
which in its turn, in view of (\ref{HessH}), is equivalent to
\[%
%TCIMACRO{\tsum \nolimits_{|J_{1}|{}=l}}%
%BeginExpansion
{\textstyle\sum\nolimits_{|J_{1}|{}=l}}
%EndExpansion
\tfrac{\partial^{2}H}{\partial p_{\alpha}^{I.i}\partial p_{\beta}^{J_{1}.j}%
}\mathfrak{s}^{\ast}(\partial_{\beta}^{J_{1}j}L-%
%TCIMACRO{\tsum \nolimits_{|K|{}\leq l}}%
%BeginExpansion
{\textstyle\sum\nolimits_{|K|{}\leq l}}
%EndExpansion
\delta_{Kk}^{J_{1}j}p_{\beta}^{K.k})=0,\quad|I|{}=l
\]
Put $E=\sum_{|K|{}\leq l}p_{\alpha}^{K.k}u_{Kk}^{\alpha}-L$. Then
\begin{align*}
\mathfrak{s}^{\ast}(E)  &  =\mathfrak{s}^{\ast}(%
%TCIMACRO{\tsum \nolimits_{|K|{}\leq l}}%
%BeginExpansion
{\textstyle\sum\nolimits_{|K|{}\leq l}}
%EndExpansion
p_{\alpha}^{K.k}u_{Kk}^{\alpha}-L)\\
&  =%
%TCIMACRO{\tsum \nolimits_{|K|{}\leq l}}%
%BeginExpansion
{\textstyle\sum\nolimits_{|K|{}\leq l}}
%EndExpansion
p_{\alpha}^{K.k}\tfrac{\partial H}{\partial p_{\alpha}^{K.k}}-\mathbb{F}\!%
\mathscr{H}^{\ast}(L)\\
&  =H.
\end{align*}
Thus
\begin{align*}%
%TCIMACRO{\tsum \nolimits_{|J_{1}|{}=l}}%
%BeginExpansion
{\textstyle\sum\nolimits_{|J_{1}|{}=l}}
%EndExpansion
\tfrac{\partial^{2}H}{\partial p_{\alpha}^{I.i}\partial p_{\beta}^{J_{1}.j}%
}\mathfrak{s}^{\ast}(\partial_{\beta}^{J_{1}i}L-%
%TCIMACRO{\tsum \nolimits_{|K|{}\leq l}}%
%BeginExpansion
{\textstyle\sum\nolimits_{|K|{}\leq l}}
%EndExpansion
\delta_{Kk}^{J_{1}i}p_{\beta}^{K.k})  &  =%
%TCIMACRO{\tsum \nolimits_{|J_{1}|{}=l}}%
%BeginExpansion
{\textstyle\sum\nolimits_{|J_{1}|{}=l}}
%EndExpansion
\tfrac{\partial^{2}H}{\partial p_{\alpha}^{I.i}\partial p_{\beta}^{J_{1}.j}%
}\mathfrak{s}^{\ast}(\partial_{\beta}^{J_{1}j}E)\\
&  =%
%TCIMACRO{\tsum \nolimits_{|J_{1}|{}\leq l}}%
%BeginExpansion
{\textstyle\sum\nolimits_{|J_{1}|{}\leq l}}
%EndExpansion
\tfrac{\partial^{2}H}{\partial p_{\alpha}^{I.i}\partial p_{\beta}^{J_{1}.j}%
}\mathfrak{s}^{\ast}(\partial_{\alpha}^{J_{1}j}E)\\
&  =[\tfrac{\partial}{\partial p_{\alpha}^{I.i}}\mathfrak{s}^{\ast
}(E)-\mathfrak{s}^{\ast}(\tfrac{\partial}{\partial p_{\alpha}^{I.i}}E)]\\
&  =[\tfrac{\partial H}{\partial p_{\alpha}^{I.i}}-\mathfrak{s}^{\ast}%
(u_{Ii}^{\alpha})]\\
&  =0.
\end{align*}
\end{proof} 

\quad

\end{document}